\documentclass[preprint]{article}

\usepackage{lineno,hyperref}
\usepackage{color}
\usepackage{amsfonts}
\usepackage{amsmath}
\usepackage{amssymb}
\usepackage{amsthm}
\usepackage{subeqnarray}
\usepackage{lscape,graphicx,url}
\usepackage{anysize}
\usepackage{proof}
\usepackage{fancyhdr}
\usepackage{latexsym}
\usepackage{array}
\usepackage{multirow}
\usepackage{textcomp}
\usepackage{optprog}
\usepackage{mathtools}
\usepackage{proof}
\usepackage[margin=2cm]{geometry}
\usepackage{setspace}
\usepackage{natbib}

\usepackage[none]{hyphenat}
\usepackage{proof}
\usepackage{paralist}
\usepackage[dvipsnames]{xcolor}
\usepackage{placeins}
\usepackage{enumitem}
\usepackage{csquotes}
\usepackage{comment}
\usepackage{bm}
\usepackage{array}
\usepackage{rotating}

\usepackage{nicematrix}

\usepackage{etoolbox}
\patchcmd{\pprintMaketitle}
  {\fi\hrule}
  {\fi\ifvoid\extrainfobox\else\unvbox\extrainfobox\par\vskip10pt\fi\hrule}
  {}{}

\newsavebox\extrainfobox

\makeatletter
\def\thickhline{%
  \noalign{\ifnum0=`}\fi\hrule \@height \thickarrayrulewidth \futurelet
   \reserved@a\@xthickhline}
\def\@xthickhline{\ifx\reserved@a\thickhline
               \vskip\doublerulesep
               \vskip-\thickarrayrulewidth
             \fi
      \ifnum0=`{\fi}}
\makeatother

\newlength{\thickarrayrulewidth}
\newcolumntype{H}{>{\setbox0=\hbox\bgroup}c<{\egroup}@{}}

\newcommand{\GG}[1]{}

\title{A two stage approach for order and rack allocation with order backlog in a mobile rack environment}

\author{Cristiano Arbex Valle$^1$ \and John E Beasley$^2$}

\date{}

\begin{document}

\maketitle

\begin{center} 
{\footnotesize

$^1$Departamento de Ci\^{e}ncia da Computa\c{c}\~{a}o, \\
 Universidade Federal de Minas Gerais, \\
 Belo Horizonte, MG 31270-010, Brasil \\
 arbex@dcc.ufmg.br \\ \vspace{0.3cm}
$^2$Brunel University \\ Mathematical Sciences, UK \\ john.beasley@brunel.ac.uk \\
}
\end{center}

\begin{abstract}
In this paper we investigate a  problem associated with operating a robotic mobile fulfilment system (RMFS). This is the problem of allocating orders and mobile storage racks to pickers. 

We present a two-stage formulation of the problem. In our two-stage approach we, in the first-stage, deal with  the orders which must be definitely fulfilled (picked), where the racks chosen to fulfil these 
first-stage orders are chosen so as to (collectively) contain sufficient product to satisfy all orders. In the second-stage we restrict attention to those racks chosen in the first-stage solution in terms of allocating second-stage orders.  

We present three different strategies for first-stage order selection; one of these strategies minimises the requirement to make decisions as to the rack sequence (i.e.~the sequence in which racks are presented to each
picker).

We present a heuristic procedure to reduce the number of racks that need to be considered.
Extensive computational results are presented for test problems that are made publicly available;  including test problems that are significantly larger than previous problems considered in the literature.
\end{abstract}


\section{Introduction}

In a facility operating as a Robotic Mobile Fulfilment System (RMFS)  mobile robots bring moveable racks of shelves containing inventory to static pickers, so that these pickers can pick the items needed for customer orders. Such a facility is an example of a
parts-to-picker system. This contrasts with the more traditional 
picker-to-parts systems, in which pickers walk/ride through the facility collecting requested items. 
In operating a RMFS decisions need to be made as to the batching of  orders (i.e.~the set of orders to be assigned to each picker), as well as which racks should be allocated (moved) to which picker. 

Robotic mobile fulfilment systems  are increasingly common  in the B2C (business to consumer) market. Moreover   their operation is an increasing topic of research within the academic community.
In this section we will not repeat material relating to the operation of a RMFS most likely familiar to the majority of readers. Rather we give here just a short description of the operation of a RMFS relevant to the problem we examine in this paper.
 Readers totally unfamiliar with the operation of a RMFS are referred to \cite{azadeh2019, boysen17a, boysen19}.

In a RMFS product to satisfy customer demand is stored on racks, also referred to as pods, or inventory pods. These racks are moved by small mobile robots to static pickers, who pick items from the rack to satisfy customer orders.
Racks are placed  in a queue leading past each picker and successively presented to the picker. The picker picks items from the presented rack and places them in one or more bins, each bin being associated with a single customer order (in some contexts these bins are known as crates or totes). 
Once all the items in an order have been picked the bin containing the completed order is put onto a conveyor, being carried away for further processing,  and a new empty bin is positioned in the space so created.  Typically  a 
scattered storage policy is adopted, so that the same product is stored in multiple racks in order to try and ensure that there is a rack containing required product free when 
needed~\citep{weidinger18b}.

This paper deals with a number of  related problems. Given a set of pickers which orders should be allocated to which picker?  Also which racks should be used so that the orders allocated to a picker can be fulfilled? Moreover in making these decisions some orders must be picked, but others need not  be. So we can have a backlog of orders that we carry forward. 
Hence we have a further decision. Which orders should be left unpicked and be backlogged? 
In this paper we present a two stage procedure to address these questions.

The structure of this paper is as follows. In 
Section~\ref{sec:review} we review the relevant literature on the operation of a robotic mobile fulfilment system. 
We also state what we believe to be the contribution of this paper to the literature.
In Section~\ref{sec:formjeb} we present our two stage formulation. We also discuss a one stage approach based upon combining our two stages.
We 
outline three different strategies that can be used with the formulation with regard to splitting the orders into a set of first-stage orders and a set of second-stage orders.
One of these strategies minimises the requirement to make decisions as to the rack sequence (i.e.~the sequence in which racks are presented to each
picker).
In Section~\ref{sec:rackred} we present our procedure to reduce the number of racks that need to be considered. In 
In Section~\ref{sec:results} we present our computational results for the test problems we examined. Finally in 
Section~\ref{sec:conclusions} we present our conclusions.

\section{Literature review}
\label{sec:review}

In this section we first review recent  literature on the operation of a robotic mobile fulfilment system. 
We also state what we believe to be the contribution of this paper to the literature.

\subsection{Recent literature}

\cite{azadeh2019, barros21,  boysen19, jaghbeer20, valle2021, zhen2022} have recently presented work discussing the literature regarding warehousing systems involving automated and robotic handling. Accordingly, for space reasons, we only consider in this section work additional to that discussed there, or work especially relevant to the problem, order and rack allocation, considered in this paper.

\cite{boysen17a} considered  the problem of a single picker with a given set of orders to be picked from a given set of racks. In such  situations decisions need to be made as to how the orders are sequenced for consideration by the picker, and too how racks are to be sequenced for presentation to the picker, given that there is a constraint on the number of orders that can be processed in parallel. They presented a mixed-integer programming model for the problem based upon the use of 
time slots. Heuristic algorithms based upon simulated annealing were also given.

\cite{li17} considered the situation where there is a single picker and formulated a zero-one integer program for the problem of deciding the racks to be allocated to the picker based upon minimising total rack travel distance (whilst supplying all orders). They also presented a three stage heuristic for the problem. 

\cite{xiang18} considered the problem of batching orders together so that all orders are assigned to some batch, and too racks are assigned to batches, which they formulated as a zero-one integer program minimising the number of racks used. They presented a heuristic to generate an initial solution as well as a variable neighbourhood search heuristic. In their work batches can be regarded as equivalent to pickers.

\cite{shi2021} presented a two stage hybrid heuristic
algorithm for the problem of simultaneously allocating orders and racks to pickers so as to minimise the number of racks used.   In the first-stage dynamic
programming together with a beam search heuristic was used to find a critical rack set, the set of
 most promising racks to use.
In the second-stage 
a constructive heuristic
 and an adaptive neighbourhood search heuristic were used to find an allocation of orders and racks to pickers. 

\cite{valle2021} considered  the problem of order and rack allocation as well as the problem of how to sequence racks for presentation to a picker using integer programs. 
They first decided the order and rack allocation and then, given the orders and racks allocated to each picker,  sequenced the racks. 
Two heuristics for order and rack allocation were given.

\cite{xie2021} considered the problem of assigning orders and racks to pickers. 
They presented a formulation of the problem as a zero-one integer program that minimises a weighted sum of the number of racks used and unused picker capacity. They do not impose a constraint which requires orders to be picked, rather the orders to be picked arise as a result of the term in the minimisation objective associated with unused picker capacity, so order
backlogs occur. They extended the formulation to account for the splitting of orders, i.e.~an order could be split so that products for the order were picked by two or more pickers. 
They presented proofs that all of the problems considered are NP-hard.
They also amended their approach to act as a heuristic for a real-world problem.

\cite{yang2021a} considered the problem as to how, for a single picker and a given set of orders, to sequence the processing of the orders  as
well as how to decide the allocation and arriving sequence of racks. 
They formulated the problem as a mixed-integer
 program and presented  a two stage solution procedure consisting of an initial greedy heuristic stage followed by an improvement stage.  

\cite{ouzidan20} considered the problem of a single picker and how to sequence the orders dealt with by that picker, as well as the racks to be brought to the picker. 
They presented a zero-one integer program for the problem, minimising the number of racks brought to the picker. They also presented a constructive heuristic for the problem which is used to generate a starting solution for a general variable neighbourhood search heuristic. 

\cite{teck2022} presented a number of mathematical programming formulations for the problems, including order and rack allocation,  encountered within a RMFS.
\cite{wang2022} considered the problem of assigning orders and racks to pickers with rack scheduling. They gave an integer programming formulation of the problem together with a two stage heuristic, assigning orders to pickers in the first-stage and then in the second-stage sequencing orders and racks. They make use of grouping by order similarity, greedy procedures and variable neighbourhood search. 

\cite{zhuang2022}  considered the problem of assigning orders and racks to pickers. 
In their paper they considered workload balancing between pickers and rack conflicts.
A rack conflict occurs when two or more pickers both require the same rack simultaneously. 
They presented a mixed-integer programming model for the problem. They also presented an adaptive large neighbourhood search heuristic that included a simulated annealing phase.

\cite{allgor2023} described the development of a  redesigned  order picking process in Amazon fulfilment centres. In their work 
not all orders need to be picked in the current \emph{\textbf{\enquote{pick window}}}, the period of time over which chosen orders will be picked. Any unpicked 
orders are added to an order backlog and hence are available to be picked in a  future pick window (along with any newly arrived orders). Their work is 
highly focused on dealing with the particular situation as encountered within Amazon fulfilment centres, e.g.~division of a fulfilment centre into multiple separate storage zones.

\cite{justkowiak2023} considered the problem of a single picker with a given set of orders to be picked from a given set of racks. They presented a mixed-integer formulation of the problem 
of deciding the order and rack sequence as well as a three-stage heuristic procedure.

In order to help identify gaps in the research literature  Table~\ref {tablep} classifies the papers discussed above.  In this table an entry of \enquote{\textbf{\emph{no}}} for multiple pickers means the paper only considers a single picker. An entry of \enquote{\textbf{\emph{no}}} for inventory constraints means the paper does not explicitly consider the number of units of each product stored on a rack, rather it assumes that if a rack contains a product then it has sufficient of that product to satisfy all orders requiring that product.  It is clear that virtually all papers use integer programming
(IP), or mixed-integer programming
(MIP),  together with heuristics.

In Table~\ref {tablep} the majority of the papers reviewed have no explicit consideration of order backlog, so all orders must be picked within the current pick window whilst 
respecting picker capacity and inventory constraints (if any). Only two papers deal explicitly with order backlog, where some orders  are left unpicked for consideration 
in  the next pick window. \emph{\textbf{It is clear from that table that there is a gap in the literature with regard to papers addressing multiple pickers with inventory constraints and order backlog. This is precisely the situation considered in this paper}}.

\begin{table}[!htb]
\centering
\renewcommand{\tabcolsep}{1mm} 
\begin{tabular}{|c|c|c|c|c|}
\hline
 & Paper  &  Multiple   & Inventory    & Solution approach \\
&  &    pickers?  &   constraints? & \\
\hline
 Order backlog & \cite{boysen17a}  & no & no    & MIP, heuristic \\
not considered & \cite{justkowiak2023}  & &  & MIP, heuristic \\
&\cite{ouzidan20} & &  &  IP, heuristic  \\
 &\cite{yang2021a}   &  &  &  MIP, heuristic \\
\cline{2-3}

& \cite{teck2022} & yes & &  MIP \\
& \cite{wang2022}   &  &   &  IP, heuristic \\
& \cite{xiang18}  &  &  &  IP, heuristic \\
& \cite{zhuang2022}  &  &    & MIP, heuristic \\

\cline{2-4}
& \cite{li17} &  no & yes &   IP, heuristic \\
\cline{2-3}
& \cite{shi2021}  & yes & &  Dynamic programming, heuristic \\
& \cite{valle2021} & &  &   IP, heuristic \\
\hline
Order backlog &\cite{allgor2023} &  yes &   yes&  MIP, heuristic \\
considered & \cite{xie2021} &   & no &   MIP, heuristic \\

\hline
\end{tabular}
\caption{Literature classification}
\label{tablep}
\end{table}


\subsection{Contribution}
\label{sec:contrib}

In this paper we have chosen to focus on order and rack allocation with multiple pickers, inventory constraints and order backlog. With regard to the importance of considering order backlog then, as emphasised by~\cite{allgor2023} with regard to Amazon fulfilment centres, \emph{\textbf{allowing a backlog of orders (so leaving some orders unpicked in the current pick window) potentially leads to better use of robot/rack resources}}. For example if in the current pick window there is a single order for a product which is stored on a rack associated with no other orders then satisfying this order would entail a robot moving the rack purely so that this single order can be fulfilled. By contrast if we allow this order to be backlogged, so moved to the next pick window, we might find that in that pick window we have newly arrived orders (or racks updated with newly arrived inventory) which would enable more effective robot/rack usage.

In light of the above literature review, as summarised in Table~\ref{tablep}, we
 believe that the contribution of this paper to the literature is:
\begin{compactitem}
\item to present a two stage optimisation based approach for simultaneously allocating both orders and racks to  multiple pickers with inventory constraints and order backlog

\item to explore use of our two stage approach with three different strategies for first/second-stage order selection; one of these strategies minimises the requirement to make decisions as to the rack sequence (i.e.~the sequence in which racks are presented to each picker)

\item to present a heuristic procedure to reduce the number of racks that need be considered; this procedure is independent of our two stage approach and hence could easily be incorporated into any other solution approach (heuristic or optimal) for order and rack allocation


\item to present extensive computational results for test problems involving multiple pickers and inventory constraints that are made publicly available for use by other researchers; including test problems that are significantly larger than previous problems considered in the literature

\end{compactitem}
With respect to the differences between the work presented in this paper and our previous work~\cite{valle2021} some elements of our first-stage formulation, Equations~(\ref{eq1})-(\ref{eq6}) below, correspond to elements seen in the problem formulation in~\cite{valle2021}. However virtually all other elements in this paper, specifically the second-stage formulation, Section~\ref{contrib1}; the combined approach, Section~\ref{jebcomb}; use of $\gamma_i$ and $\delta_o$, Section~\ref{contrib2}; choice of $F$ and $S$, Section~\ref{sec:threejeb};
rack reduction, Section~\ref{sec:rackred}; plus the material in~\ref{jebappa} and~\ref{jebappb} together with the associated computational results
in Section~\ref{sec:results};
 were not seen previously in~\cite{valle2021}.

\section{First-stage and second-stage formulations}
\label{sec:formjeb}
Suppose that there are in total  $N$ products within a facility  that can be ordered by customers.  In the current pick window we have $O$ orders to be allocated to pickers where each order $o$ requires $q_{io}$ units of product $i$. 
In a practical setting it is very likely that the set of orders currently being considered does not involve all of the $N$ products, so let the set of active products $I$ be defined using $I = [i~|~ \sum_{o=1}^O q_{io} \geq 1~ i =1,\ldots,N]$.
Let $Q$ be the set of orders comprising just a single unit of one product, so $Q = [o~|~ \sum_{i \in I} q_{io} = 1~ o =1,\ldots,O]$.
There are $P$ pickers, where picker $p$   has a capacity for at most $C_p$ orders.
Without loss of generality assume that the pickers are indexed in decreasing $C_p$ order (so $C_p \geq C_{p+1}~p=1,\ldots,(P-1)$). 

We have $R$ mobile racks which can be allocated (moved)  to pickers.
Rack $r$ contains $s_{ir}$ units of product $i$. Let $\Gamma(i)$ be the set of racks that contain one or more units of product $i$, so $\Gamma(i) = [r~|~s_{ir} \geq 1~ r=1,\ldots,R]$.  
Let $\Delta(o) = [r~|~r \in \Gamma(i)~q_{io} \geq 1~i \in I]$ be the set of racks that involve one (or more) products in order $o$.
Let $\gamma_i$ be the minimum number of racks that are needed to supply all orders involving product $i$. 
Let $\delta_o$ be the minimum number of racks that are needed to supply order $o$.
We indicate later below how to find  values for $\gamma_i$ and $\delta_o$.
We assume without loss of generality that any orders that cannot be satisfied at all, 
even by using the entire set of racks, have been identified and eliminated from the problem.

In our approach we partition the set of orders into two disjoint sets $F$ and $S$ where $F \cap S= \emptyset$ and $F \cup S = [o~|~o=1,\ldots,O]$. 
$F$ is the set of orders which are of most importance and must be fulfilled (picked), so allocated to pickers whilst respecting picker capacity constraints.   
$S$ is the set of  orders which are of lesser importance and may potentially be left unpicked (so constituting the backlog for the next pick window).
For example (in practice) orders may be  in $S$ if the latest time by which their pick should be scheduled such that the customer receives their order on time has not yet been reached.
Our approach involves:
\begin{compactitem}
\item  a first-stage formulation allocating both orders $o \in F$ and racks to pickers whilst  maintaining inventory levels sufficient for orders $o \in S$
\item  a second-stage formulation allocating (if possible) orders $o \in S$ to pickers, but restricting attention to just the set of racks assigned to pickers at the first-stage
\end{compactitem}
The logic here is that the set $F$ of orders must be picked in the current pick window, and in the first-stage racks are chosen so as to do this. 
But given that these racks are allocated to pickers it is clearly of benefit if the inventory on these racks is used to fulfil orders in $S$, and this is what the second-stage addresses.
\emph{\textbf{Any unpicked orders left after the first and second stage constitute the order backlog for the consideration in the next pick window}}. This backlog, together with any newly arrived orders, is  partitioned  into two (new) disjoint sets $F$ and $S$ and the process repeats over a succession of pick windows.

\subsection{First-stage formulation}
Our first-stage formulation for allocating orders and racks to pickers involves the following zero-one decision variables:
\begin{compactitem}
\item $x_{op}=1$ if order $o$ is allocated to picker $p$, zero otherwise
\item $u_r=1$ if  rack $r$ is used, so allocated to some picker, zero otherwise
\item $y_{rp}=1$ if  rack $r$ is allocated to picker $p$, zero otherwise
\end{compactitem}


\noindent Our  formulation is:
\begin{optprog}
\optaction[]{min} &  \objective{M \sum_{r=1}^R  u_r -
\sum_{i \in I} (
(\sum_{r=1}^R s_{ir}u_r - \sum_{o=1}^O q_{io} )/\sum_{o=1}^O q_{io})
\label{eq1} }\\

subject to: & \sum_{o \in F} x_{op} & \leq & C_p & \forall p \in \{1,\ldots,P\}
\label{eq2} \\

& \sum_{p=1}^P x_{op} & = & 1 & \forall o \in F
\label{eq2a} \\

& \sum_{p=1}^P y_{rp} & = & u_r & \forall r \in \{1,\ldots,R\}
\label{eq3} \\

& \sum_{r=1}^R s_{ir}y_{rp} & \geq & \sum_{o \in F}q_{io}x_{op} & 
\forall i \in I;  
p \in \{1,\ldots,P\}:~\sum_{o \in F} q_{io} \geq 1
\label{eq4} \\

& \sum_{r=1}^R s_{ir}u_r& \geq & \sum_{o=1}^O q_{io} & 
\forall i \in I
 \label{eq4aa} \\

& \sum_{r=1}^R y_{rp} & \geq & 1 &  \forall p \in \{1,\ldots,P\}
 \label{eq4u} \\

& \sum_{r \in \Gamma(i)} u_r & \geq & \gamma_i & 
\forall i \in I
 \label{eq4a} \\

& \sum_{r \in \Delta(o)} u_r & \geq & \delta_o & 
\forall o \in F
 \label{eq4b} \\

& u_r & = & 0 & \forall r \in \{1,\ldots,R\}:~\sum_{o=1}^O \sum_{i \in I} s_{ir}q_{io}=0
\label{eq7}  \\

& x_{op} & \in & \{0,1\} & \forall o \in F;  p \in \{1,\ldots,P\}
\label{eq5} \\
& u_r & \in & \{0,1\} &\forall r \in \{1,\ldots,R\}
\label{eq5a} \\
& y_{rp} & \in & \{0,1\} & \forall r \in \{1,\ldots,R\}; p \in \{1,\ldots,P\}
\label{eq6} 
\end{optprog}

Equation~(\ref{eq2}) ensures that we allocate orders $o \in F$ to each picker so as to respect  picker capacity. 
Equation~(\ref{eq2a}) ensures that each order $o \in F$ is allocated to one picker. Note here that in this first-stage formulation we are not allocating orders $o \in S$ to pickers. Equation~(\ref{eq3}) ensures that each rack is either allocated to one picker, or not used at all.
Equation~(\ref{eq4}) is the product supply constraint and ensures that, for each product $i$ and each picker $p$, the number of units of that product available from the racks assigned to that picker (so  $ \sum_{r=1}^R s_{ir}y_{rp}$) is sufficient to meet the required number of units of product $i$ at picker $p$ given the orders in $F$ allocated to the picker (so  $ \sum_{o \in F}q_{io}x_{op}$). Hence Equation~(\ref{eq4})  ensures that appropriate racks are assigned to picker $p$ so as to enable  the picking of all of the products associated with the orders $o \in F$ assigned to the picker. 

Equation~(\ref{eq4aa}) ensures that, over the entire set of racks chosen, the total number of units of product $i$ available in those racks (so $\sum_{r=1}^R s_{ir}u_r$)
is sufficient to satisfy all orders (where the orders require in total $\sum_{o=1}^O q_{io}$ units of product $i$). This constraint links the racks chosen at the first-stage with the allocation of orders $o \in S$ to these chosen racks at the 
second-stage in terms of having sufficient of each product available (after allocation of orders $o \in F$) for orders $o \in S$ at the 
second-stage  (as in  the second-stage formulation presented later below).

Equation~(\ref{eq4u}) means that we choose at least one rack for each picker. The assumption behind this constraint is that since we have $P$ pickers available we wish to make use of them.
Equation~(\ref{eq4a}) ensures that for each product $i$ we use at least the minimum number ($\gamma_i$) of racks needed.  Equation~(\ref{eq4b}) ensures that for each order $o \in F$, which we know (from Equation~(\ref{eq2a})) must be allocated to some picker, we use at least the minimum number ($\delta_o$) of racks needed.

Given the $O$ orders then we know the products needed. If a rack contains  none of the required products then clearly it is irrelevant and will never be assigned to a picker, and so can be removed from the problem. Equation~(\ref{eq7}) ensures that any such  racks are eliminated  since 
the term $\sum_{o=1}^O \sum_{i \in I} s_{ir}q_{io}$ is zero if and only if rack $r$ does not contain any of a required product.
Equations~(\ref{eq5})-(\ref{eq6}) are the integrality constraints.

In our objective, Equation~(\ref{eq1}), $M$ is a large positive constant and the first term in that objective ensures that we minimise the total number of racks allocated to pickers. Using  an objective related to rack minimisation is very common in the literature~(\cite{boysen17a, hansen18, justkowiak2023, 
ouzidan20,
shi2021, valle2021, xiang18, xie2021}) as it contributes to less movement of racks from the storage area to the picking area, and too encourages the picking of multiple products for different orders from the same rack. 

Given that there may be alternative optimal solutions involving the same (minimum) number of racks we make use of the second term in the objective. This term is designed to maximise the number of units of products $i \in I$ associated with the racks chosen. 
In this objective the term  
$(\sum_{r=1}^R s_{ir}u_r - \sum_{o=1}^O q_{io} )$ is the total number of units of product $i$ on the racks chosen $(\sum_{r=1}^R s_{ir}u_r)$
in excess of those needed $(\sum_{o=1}^O q_{io} )$, where this excess is proportionally scaled by the total number needed $(\sum_{o=1}^O q_{io} )$. This term, with its associated negative sign in Equation~(\ref{eq1}), 
means that we maximise the number of units of products $i \in I$ associated with the racks chosen. This gives us additional flexibility in allocating orders $o \in S$ at the second-stage. 

Rearranging Equation~(\ref{eq1}), and dropping the constant term, we have that our objective is:
\begin{optprog}
\optaction[]{min} &  \objective{\sum_{r=1}^R  (M -
\sum_{i \in I} 
(s_{ir}/\sum_{o=1}^O q_{io})) u_r
\label{eq1simp} }
\end{optprog}
With this objective  a suitable value for $M$ is any value that exceeds $\sum_{r=1}^R  
\sum_{i \in I} (
s_{ir}/\sum_{o=1}^O q_{io})$.

In the first-stage formulation presented above Equation~(\ref{eq4aa}) is included to link with the second-stage and ensure we have sufficient inventory for all orders. Equations~(\ref{eq4u})-(\ref{eq7}) are valid inequalities included to strengthen the linear programming relaxation.

\subsection{Second-stage formulation}
\label{contrib1}

Once the first-stage formulation presented above, Equations~(\ref{eq1})-(\ref{eq6}), has been solved then we will have decided the racks to be allocated to each picker. In addition, for all orders $o \in F$, the picker to which they have been allocated will also have been decided. 
In our approach  \emph{\textbf{the racks chosen at the first-stage are not changed when we come to the second-stage}}. 
Allocation of all orders at the second-stage then occurs only considering the set of racks chosen at the first-stage (rather than considering the complete set of racks).
The logic behind this is that we know from the first-stage solution that with this chosen set of racks we can feasibly allocate all orders $o \in F$ (from Equation~(\ref{eq2a})). However, by judicious reallocation of orders to racks, and racks to pickers, we may be better able to deal with allocating  orders $o \in S$ to pickers, as these orders have yet to be allocated.

Equation~(\ref{eq4aa}) ensures that collectively the racks chosen in the first-stage formulation solution have sufficient inventory left (after dealing with orders $o \in F$) to deal with all orders $o \in S$. 
However, this is not sufficient to ensure that we can always feasibly allocate all orders $o \in S$ to pickers. This may, for example, be due to the fact 
that the racks allocated to a picker do not contain sufficient product to satisfy an order, even if the racks in total have sufficient product available (recall here that orders must be supplied by a single picker). Alternatively allocation of all orders $o \in S$ to pickers whilst respecting rack inventory positions may require violation of picker capacity.
To account for these possibilities we, in allocating orders $o \in S$ to pickers, aim to maximise the number of allocated orders.
Note here that use of the first-stage objective function term associated with maximising the number of units of products $i \in I$ was designed to enable more orders to be allocated to pickers at the second-stage.

Let $\Theta = [r~|~u_r=1~r=1,\ldots,R]$ be the set of racks chosen in the first-stage solution. 
Let $v_o=1$ if order $o \in S$ is allocated to some picker in the second-stage solution, zero otherwise. 
Then our second-stage formulation is:
\begin{optprog}
\optaction[]{max} &  \objective{ \sum_{o \in S}  v_o
\label{feq1} }\\

subject to: & \sum_{o=1}^O x_{op} & \leq & C_p & \forall p \in \{1,\ldots,P\}
\label{s2eq2} \\

& \sum_{p=1}^P x_{op} & = & 1 & \forall o \in F
\label{s2eq2a} \\

& \sum_{p=1}^P x_{op} & = & v_o & \forall o \in S
\label{s2eq2b} \\

& \sum_{p=1}^P y_{rp} & = & 1 & \forall r \in \Theta
\label{s2eq3} \\

& \sum_{r \in \Theta} s_{ir}y_{rp} & \geq & \sum_{o=1}^O q_{io}x_{op} & 
\forall i \in I;  p \in \{1,\ldots,P\}
\label{s2eq4} \\

& v_o & \in & \{0,1\} & \forall o \in S
\label{s2eq4a} \\ 
& x_{op} & \in & \{0,1\} & \forall o \in \{1,\ldots,O\};  p \in \{1,\ldots,P\}
\label{s2eq5} \\
& y_{rp} & \in & \{0,1\} & \forall r \in \Theta; p \in \{1,\ldots,P\}
\label{s2eq6} 
\end{optprog}

Equation~(\ref{feq1}) maximises the  number of orders $o \in S$ allocated to pickers.  
Equation~(\ref{s2eq2}) ensures that for picker $p$ the capacity $C_p$ is respected. Equation~(\ref{s2eq2a}) ensures that each order $o \in F$ is allocated to one picker and
Equation~(\ref{s2eq2b}) ensures that each order $o \in S$ is either allocated to one picker (or not allocated).
Equation~(\ref{s2eq3}) ensures that each rack $r \in \Theta$ is allocated to a picker. Equation~(\ref{s2eq4}) is the product supply constraint and ensures that, for each product $i$ and each picker $p$, the number of units of that product available from the racks assigned to that picker  is sufficient to meet the required number of units of product $i$ at picker $p$ given the orders   allocated to the picker.
Equations~(\ref{s2eq4a})-(\ref{s2eq6}) are the integrality constraints.

There are a number of extensions to these first and second stage formulations to deal with issues that may arise in practice
and these  are detailed in~\ref{jebappa}.
\subsection{Combined approach}
\label{jebcomb}
Above we have presented a two stage procedure with two related optimisation problems being solved in a sequential fashion. Our first-stage essentially minimises the number of racks chosen, the second term in Equation~(\ref{eq1}) being a tie break between solutions with equal minimum number of racks. Then in the second-stage we maximise the number of orders in $S$ we can deal with from the racks chosen in the first-stage, given that we must pick orders in $F$. It is possible to combine the two approaches with minor modifications. Consider the combined approach:
\begin{equation}
\mbox{min} ~~~ (|S|+1) \sum_{r=1}^R u_r - \sum_{o \in S} v_o \label{c1} 
\end{equation}
subject to
Equations~(\ref{eq2})-(\ref{eq6}),(\ref{s2eq2}),(\ref{s2eq2b}),(\ref{s2eq4}),(\ref{s2eq4a}), where in Equation~(\ref{s2eq4}) the set $\theta$ is now taken to be the full rack set. Here $(|S|+1)$ is a suitable weighting factor to give priority to the first term in this combined objective, the number of racks used, as in our first-stage objective (Equation~(\ref{eq1})) above. The effect of the  negative weighting on the $\sum_{o \in S} v_o$ term is to maximise the number of orders picked given the racks chosen, as in our second-stage objective (Equation~(\ref{feq1})) above.

This combined approach would need to consider $OP$ $x$ variables and $RP$ $y$ variables. However in our two stage approach we solve smaller problems; so in the first-stage a problem with $|F|P $ $x$ variables and $ RP $ $y$ variables and in the second-stage a problem with $ OP $ $x$ variables and $|\theta|P $ $y$ variables. Our two stage approach is based on the assumption that it will be computationally more effective to solve two smaller problems as compared with solving one larger (combined) problem. Computational results, presented later below, support this assumption. 
Note here that it is a simple extension to the combined approach to alter the relative  weighting given to the two summation terms in 
Equation~(\ref{c1}), so that the objective becomes minimise 
$\lambda_1 \sum_{r=1}^R u_r - \lambda_2\sum_{o \in S} v_o$, where $\lambda_1,\lambda_2 >0$ are the weights given to each of the summation terms.


\subsection{Finding $\gamma_i$ and $\delta_o$}
\label{contrib2}
Identifying a value for $\gamma_i$, the minimum number of racks that are needed to supply all orders involving product $i$,
 can be done by solving a simple 
zero-one integer program involving just a single constraint. For product $i \in I$ we have that $\gamma_i$ corresponds to the optimal value of:
\begin{optprog}
\optaction[]{min} &  \objective{\sum_{r \in \Gamma(i)}  u_r  \label{eqg1} }\\

subject to: & \sum_{r \in \Gamma(i)} s_{ir}u_r & \geq & \sum_{o=1}^O q_{io}  
\label{eqg2} \\

& u_r & \in & \{0,1\} &\forall r \in \Gamma(i)  \label{eqg3} 
\end{optprog}
Equation~(\ref{eqg1}) minimises the number of racks used whilst Equation~(\ref{eqg2}) ensures that for the product considered sufficient racks are chosen to be able to supply all orders with that product. Equation~(\ref{eqg3}) is the integrality constraint.
Equations~(\ref{eqg1})-(\ref{eqg3}) is a  simple zero-one problem.

In a similar fashion 
we have that $\delta_o$, the minimum number of racks needed to fulfil order $o \in F$, corresponds to the optimal value of:
\begin{optprog}
\optaction[]{min} &  \objective{\sum_{r \in \Delta(o)} u_r\label{jebd1} }\\

subject to: 
& \sum_{r \in \Delta(o)} s_{ir}u_r & \geq & q_{io}
& \forall i \in I:~ q_{io} \geq 1  
\label{jebd2} \\

& u_r & \in & \{0,1\} &\forall r \in \Delta(o)  \label{jebd3} 
\end{optprog}
Equation~(\ref{jebd1}) minimises the number of racks used whilst Equation~(\ref{jebd2}) ensures that for each product involved in order $o$ we have sufficient inventory on the racks used to supply that product. Equation~(\ref{jebd3}) is the integrality constraint.
Equations~(\ref{jebd1})-(\ref{jebd3}) is a simple zero-one problem.

In terms of computational effort for the optimisation problems given above we need to solve $|I|$ cases of the zero-one program Equations~(\ref{eqg1})-(\ref{eqg3})  and $|F|$ cases of zero-one program Equations~(\ref{jebd1})-(\ref{jebd3}). 
Theoretically both of these zero-one optimisation programs could require a solution time exponentially related to problem size. However  
our computational experience has been that solving these zero-one programs did not require a significant amount of time. Over all instances considered later with $O \leq 500$, the maximum time required to solve these $|I| + |F|$ optimisation problems was only 1.09 seconds,  out of a total computation time of 300 seconds. Hence effectively 
approximately~\emph{\textbf{only one-third of one percent of total computation time}}.

\subsection{Choice of $F$ and $S$}
\label{sec:threejeb}

Our two stage approach requires us to have a strategy to define the set of orders $F$ that are allocated to a picker at the first-stage, or the set of orders $S$ dealt with at the second-stage (recall here that $F \cap S= \emptyset$ and $F \cup S = [o~|~o=1,\ldots,O]$, so defining one of these two sets automatically defines the other). Clearly there can be many possible strategies for splitting the set of orders into two mutually exclusive subsets $F$ and $S$. 

In any particular practical application an individual RMFS operator would have their own distinct strategy for defining $F$ and $S$.
\emph{\textbf{However we should stress here that in the work reported this paper we are not working with any commercial RMFS partner who might have their own individual  strategy for defining $F$ or $S$ which could be investigated computationally. }} For this reason we are unable to investigate tailor-made strategies designed for use within a particular RMFS facility. 
However in this paper we do consider three different strategies for choice of the sets $F$ and $S$, each with their own underlying logic.  

Strategies other than the three considered below are possible. For example randomly assigning orders to $F$ and $S$.  However such a random allocation  lacks (in our view) a clear defining logic as to why it might be adopted. Obviously  space reasons require us to limit the number of possible strategies that we can examine and report here. The three strategies that we did examine are:
\begin{compactitem}
\item \emph{\textbf{Strategy 1}}: $S = \emptyset$.
\item \emph{\textbf{Strategy 2}}: $S = Q$, so $S$ is the set of orders comprising just a single unit of one product
\item \emph{\textbf{Strategy 3}}: $F = [o~|~ \delta(o)=1 ~ o=1,\ldots,O]$, so $F$ is the set of orders which only require one rack
\end{compactitem}
\noindent The logic underlying each of these three strategies is explained below.

\subsubsection{Strategy 1}

In Strategy 1 since $S = \emptyset$ we allocate all orders to pickers in the first-stage, in other words we pick all orders.
 As such there are no decisions to be made at the second-stage. 
We would expect this strategy to be computationally the most demanding strategy  as it involves explicit consideration of all orders $F=[o~|~o=1,\ldots,O]$ in the first-stage. 
Note here that since for this strategy we have no second-stage the second term in the objective (Equation~(\ref{eq1})) becomes irrelevant and so for this strategy we simply minimise $\sum_{r=1}^R u_r$. 

\emph{\textbf{The logic underlying Strategy 1  is that we pick all orders 
 in the first-stage}}.

\subsubsection{Strategy 2}
With regard to Strategy 2 then in our previous  work~\cite{valle2021} we showed that when formulating the problem of  rack sequencing at each picker  advantage was gained  by explicit and separate consideration of orders comprising a single unit of one product  (subject to certain conditions being satisfied). Distinguishing such orders seems appropriate for environments, such as Amazon~\citep{weidinger18c}, where the vast number of orders are for only one or two items. 

Recall here that, as discussed above, order and rack allocation is followed (for each picker individually) by the sequencing of racks so that the picker can fulfil the orders assigned to them. If an order involves just a single unit of one product then we automatically know that it can be satisfied by just a single rack. All other orders, including any requiring two or more units of a single product, \emph{\textbf{may}} require two or more racks. 
For example an order for two units of the same product may take one unit from two different racks. An order for one unit of each of  two different products may also take  one unit for these two products from two different racks. 

\emph{\textbf{The logic underlying Strategy 2 is that we deal with the \enquote{difficult to pick} orders potentially needing  two or more racks in the first-stage, leaving the \enquote{easy to pick} orders, so those that just require a single unit of a single product, to the second-stage}}.

\subsubsection{Strategy 3}

In Strategy 3 $F$ is the set of orders which can be satisfied (fully supplied) by a single rack (i.e.~the set of orders $o$ with $\delta_o=1$). In this strategy we first deal with all orders that only require (as a minimum) one rack. The remaining orders (which by definition require two or more racks) will be  satisfied (if possible) at the second-stage.  
We would note here $Q \subseteq F$, since clearly the set $Q$ of 
single unit, single product, orders can be fully supplied  by a single rack. 
  However $F$ can also involve other orders that can be fully supplied by just a single rack (e.g.~an order for two different products where sufficient of both products are stored on the same rack).

In general terms this strategy aims to streamline picker rack sequencing. Recall here, from the description of the problem given above, that the racks allocated to a picker must be sequenced such that all orders allocated to a picker are picked using the bins available for assembling orders at a picker. If it is known that all of the products (and their appropriate quantities) for an order can be fully supplied (picked) from a single rack then this helps simplify the sequencing of picking operations. 

For this strategy we do need to introduce additional constraints into the formulation, to ensure that all
 orders which can be satisfied (fully supplied) by a single rack are indeed supplied from just a single rack,
 and these are detailed in~\ref{jebappb} together with 
 a number of other constraints that we can introduce into the formulation.

\emph{\textbf{The logic underlying Strategy 3  is that all orders that can be satisfied by a single rack are picked using just a single rack in the first-stage, thereby minimising  the requirement to make decisions as to the rack sequence, 
leaving orders that involve two or more racks to the second-stage}}.


\

\section{Rack reduction}
\label{sec:rackred}

In our first-stage formulation we potentially consider all racks (although Equation~(\ref{eq7}) eliminates from consideration any racks that cannot contribute any product items for the orders considered). In order to 
 improve computational performance we can reduce the number of racks that we consider in the first-stage, selecting for consideration only a small set of racks. The procedure we adopted for this is given below.

In terms of rack selection for the first-stage formulation (Equations~(\ref{eq1})-(\ref{eq6})) we aim to include sufficient racks to enable us to fulfil orders $o \in F$ (Equation~(\ref{eq2a})). In addition we aim to include racks to enable the fulfilment of all orders (Equation~(\ref{eq4aa})).  

In order to present our rack selection procedure for the order and rack allocation problem let $\Phi$ be the set of possible racks for selection, where initially $\Phi \leftarrow [1,\ldots,R]$.
Let $\Lambda$  be an ordered (sequenced) set of customer orders, with the first orders in this list being customer orders $o \in F$ (in descending 
$\sum_{i \in I} (q_{io}/\sum_{r \in \Phi} s_{ir})$ 
order, ties broken arbitrarily) and the remaining orders in this list being customer orders $o \in S$ (also in descending 
$\sum_{i \in I} (q_{io}/\sum_{r \in \Phi} s_{ir})$ 
order, ties broken arbitrarily). Here the term $(q_{io}/\sum_{r \in \Phi} s_{ir})$ represents for order $o$ the fraction that $q_{io}$ is of the total inventory available for product $i$ on the racks considered, $\sum_{r \in \Phi} s_{ir}$. 

Define $T_r$ as the total number of orders that it is possible to (potentially) fulfil from rack $r$, but allowing for fractional fulfilment of an order. So:
\begin{optprog}
& & T_r & = & \sum_{o \in \Lambda} ~\Bigl\{ ~(\sum_{i \in I} \mbox{min}[s_{ir},q_{io}])/
 (\sum_{i \in I} q_{io} ) ~\Bigr\} & \forall r \in \Phi
\label{jebh1} 
\end{optprog}
 In Equation~(\ref{jebh1}) the term
$(\sum_{i \in I} \mbox{min}[s_{ir},q_{io}])/
 ( \sum_{i \in I} q_{io} )$ 
can be regarded as the proportion of order $o$ that can be supplied by rack $r$ given the set $\Lambda$ of orders we have to consider.

Let $\pi_{ip}$ be the number of units of product $i$  available at picker $p$. Set $\pi_{ip} \leftarrow 0, ~ \forall i \in I,~p=1,\ldots,P$ as initially no racks are allocated to any picker. Let $\Upsilon$ be the set of racks that we select, where initially $\Upsilon \leftarrow \emptyset$.
Then our procedure for selecting racks is:
\begin{enumerate}[label=(\alph*), nosep]

\item Repeat until all racks have been considered:
\newline
For the rack $r$ with the maximum value of $T_r$ (ties broken arbitrarily):

\begin{enumerate}[label=(\roman*), nosep]
\item Allocate the rack $r$ currently being considered to the picker $p$ with maximum $C_p$ value (ties broken arbitrarily). Set $\pi_{ip} \leftarrow \pi_{ip} + s_{ir} ~\forall i \in I$.

\item Consider the customer orders $o \in \Lambda$ in turn and allocate each such order to the current picker $p$ if the order can be completely fulfilled using the racks now allocated to picker $p$, i.e.~if $C_p \geq 1$ and $q_{io} \leq \pi_{ip}~\forall i \in I: q_{io} \geq 1$. If the order is allocated to picker $p$ then set  $C_p \leftarrow C_p - 1$;  $\pi_{ip} \leftarrow \pi_{ip} - q_{io} ~\forall i \in I: q_{io} \geq 1$; 
$ \Lambda \leftarrow \Lambda \setminus [o]$; $\Upsilon \leftarrow \Upsilon \cup [r]$.

\item If after considering all customer orders $o \in \Lambda$ in turn no order has been allocated to picker $p$ then it is clear that rack $r$ has made no contribution to the fulfilment of orders so we remove it from picker $p$, i.e.~set $\pi_{ip} \leftarrow \pi_{ip} - s_{ir} ~\forall i \in I$.
\end{enumerate}

\item Consider all products $ i \in I$ in turn and if there exists any product $i$ for which there is insufficient of the product on the racks selected to supply total demand (i.e.~if $\sum_{r \in \Upsilon} s_{ir} < \sum_{o=1}^O q_{io}$)
then consider the racks in  descending $T_r$ order and add racks $r \notin \Upsilon$ for which $s_{ir} \geq 1$ to $\Upsilon$ until there is sufficient available. 

\item If we are considering Strategy 3 then consider all orders  $o \in F \setminus Q$ and if $\Omega(o) \cap \Upsilon = \emptyset$ add the rack $r \in \Omega(o)$ with the maximum $T_r$ value 
(ties broken arbitrarily) 
to $\Upsilon$.

\end{enumerate}

In the above procedure we consider racks in decreasing order of the number of (potential) customer orders they can deal with and allocate each such rack to the picker that has maximum capacity to deal with customer orders.
We then allocate any customer orders that can be completely fulfilled from the set of racks available at the picker to which the newly allocated rack has been assigned. 
We then add racks to ensure that we have sufficient of each product on the racks selected to supply all orders. 
In the special case of Strategy 3 we add a rack for each order than does not has a rack available to fully supply it.

We would stress here that the rack reduction procedure we have presented above could easily be incorporated into any other solution approach (heuristic or optimal) for order and rack allocation, e.g.~by applying it with $F = [o~|~ o=1,\ldots,O]$, so $F$ is the entire set of orders.

To incorporate our rack reduction procedure into our 
two stage approach it is clear that we need in the first-stage formulation to restrict attention to the set  $\Upsilon$ of racks selected. This can be easily done by adding the constraint
$u_r = 0 ~\forall r \notin \Upsilon$
 to the first-stage formulation. However excluding from consideration all racks $ r \notin \Upsilon$ may lead to infeasibility, and so we do include all racks, but penalise in the objective function use of any rack $ r \notin \Upsilon$.
For the very large problems considered below (with $O \geq 3200$) we found it computationally beneficial to amend step (a) of the above rack reduction procedure so that if, when considering rack $r$, we allocate an order $o \in F$ to picker $p$ we also set $x_{op}=1$ and $y_{rp}=1$. In addition we carry forward to the second-stage variable values from the first-stage solution.

\section{Computational results}
\label{sec:results}

In this section we present results for the 
solution of our two stage approach to order and rack allocation.

\subsection{Test problem instances}

\sloppy As discussed above our two stage approach (using Strategies 2 and 3) does not require all orders to be picked, rather some can be left unpicked. We would expect that if total picker capacity is greater than the total number of orders then it would be easier to fulfil all orders. Hence in the computational results presented below we only consider instances in which total picker capacity is equal to the total number of orders, i.e.~$\sum_{p=1}^P C_p = O$. The test instances examined, for $O=50,100,150,200,500$, were generated in the manner described in full detail  in~\cite{valle2021}, in particular with
 the average order comprising just 1.6 items~(\cite{boysen19, boysen19b}) and the vast number of orders containing only one or two items~(\cite{weidinger18c}).
These problem instances are publicly available for use by other workers and can be found at
\href{https://www.dcc.ufmg.br/~arbex/mobileRacks.html}{https://www.dcc.ufmg.br/$\sim$arbex/mobileRacks.html}.

In the computational results presented below we used an Intel Core i7-7500U @ 2.70GHz with 8GB of RAM and Linux as the operating system. The code was written in C++ and 
Cplex 12.10~\citep{cplex1210} was used as the 
solver.
A maximum time limit of 5 CPU minutes (300 seconds) was imposed unless otherwise stated. This is the same time limit as in our previous work~\cite{valle2021}. As stated there the logic for imposing a 5 minute time limit was that we might reasonably expect that in a B2C environment decisions as to order and rack allocation 
for a pick window
have to be made relatively quickly.

\subsection{Order and rack allocation: no rack reduction}

Table~\ref{nored} shows the results for $O=50,100,150,200,500$ when we apply our varying strategies for order and rack allocation with no rack reduction. We also give here results for directly applying the approach given in our previous work~\cite{valle2021} to the same set of test problems.

In the computational results tables shown below {\it\textbf{T(s)}} denotes the total computation time in seconds and 
{\it\textbf{UB}} is the best solution (upper bound) obtained at the end of the search, either when the instance was solved to proven optimality or when the time limit was reached. {\it\textbf{GAP}} is defined as $100 (\text{UB} - \text{LB})/\text{UB}$ where 
LB is the best lower bound obtained at the end of the search, either when the instance was solved to proven optimality or when the time limit was reached. If the time limit is reached then this is indicated by {\it\textbf{TL}} in the tables below.

For some of the  problems reported below no upper bound (feasible solution) was found before the problem terminated at time limit, and consequently for these problems no GAP value can be reported. Problems where this occurs
are signified    by 
{\it\textbf{NF}} for No Feasible solution.

In the results for Strategies 1, 2 and 3 {\bm{$|\Theta|$}} is the number of racks chosen in the first-stage solution and {\bm{$|S|$}} is the number of orders that are considered for picking in the second-stage. {\it\textbf{BLG}} is the {\it\textbf{BackLoG}}, so the number of 
unpicked orders, i.e.~orders not picked at the second-stage that must be carried over to the next pick window.

To illustrate the results consider the results for $O=50$ orders; $N=500$ products; $R=75$ racks;  $P=10$ pickers; each with a capacity $C_p$ of 5 orders (so the last line in the upper body of Table~\ref{nored}, related to $O = 50$):
\begin{compactitem}

\item The approach given in our previous work (\cite{valle2021}), when run on the same hardware with the same version of Cplex as used in this new work, requires 19.3 seconds, has a zero GAP 
(so solves the problem to proven optimality) with the optimal solution requiring 23 racks.  This means that 23 racks was the minimum number of racks needed to deal with the total of $\sum_{p=1}^P C_p=50$ orders allocated to the pickers.

\item Strategy 1, which requires all orders to be picked, solves the problem to the same proven optimal solution, requiring 70.3 seconds. For Strategy 1 the second-stage is irrelevant as all orders are picked at the first-stage. 

\item Strategy 2, which leaves orders requiring just a single unit of one product to the second-stage, also solves the problem to proven optimality. The first-stage required 17.5 seconds and the second-stage 0.3 seconds. 21 racks were needed, with the number of second-stage orders being 28. 7 of these 28 orders were left unpicked. This is because in the second-stage we restrict ourselves to only using the racks associated with the first-stage solution. 

\item Strategy 3  requires orders  which can be supplied from a single rack  to be picked at the first-stage. 
 The first-stage required 0.1 seconds and the second-stage 0.4 seconds. 21 racks were needed, with the number of second-stage orders being 14. 3 of these 14 orders were left unpicked. As for Strategy 2 these orders are left unpicked because in the second-stage we restrict ourselves to only using the racks associated with the first-stage solution. 

\end{compactitem}

\noindent Considering Table~\ref{nored} then, as we would expect, as problem size (in terms of the number of orders) increases we have an increasing number of instances terminating due to reaching the time limit. 

Both~\cite{valle2021} and Strategy 1 require all orders to be fulfilled in the first-stage. Referring to the averages seen at the foot of each of these tables then for $O \leq 150$  the approach 
of~\cite{valle2021} appears to perform better than Strategy 1, given the time limit of 300 seconds imposed. It requires the same (or less) time whilst returning a solution involving the same, or fewer, racks. However that performance advantage disappears for $O=200$ and $O=500$. For these problems (all of which went to time limit) we have that Strategy 1 required (on average) fewer racks, indeed significantly fewer racks when $O=500$.

Strategy 2, which leaves orders requiring only a  single unit of one product to the second-stage, appears to be increasingly effective as the number of orders increases. Indeed for $O=500$ it (on average) gave solutions requiring far fewer racks than either~\cite{valle2021} or Strategy 1, and  over all ten cases regarding $O = 500$ in Table~\ref{nored} left no backlog.

Strategy 3, which in the first-stage deals with the set of orders that can be satisfied (fully supplied) by a single rack, on average requires more racks than Strategy 2. However it is important to be clear here this result, requiring more racks, is an inherent feature of Strategy 3. Recall that Strategy 3 chooses racks  such that all first-stage orders are fully satisfied by
just  a single rack. Clearly having sufficient racks such that all first-stage orders can be fully satisfied by just one of the chosen racks will be expected to require more racks than if we are allowed (as in Strategies 1, 2 and~\cite{valle2021}) to satisfy orders using products picked from a number of different  racks.

\begin{table}[!htpb]
\centering
{\tiny
\renewcommand{\tabcolsep}{0.8mm} 
\renewcommand{\arraystretch}{1.8} 
\begin{tabular}{|c|c|c|c|c|rrr|rrr|rrrrrrr|rrrrrrr|}
\hline
 \multirow[c]{3}{*}{$O$} &  \multirow[c]{3}{*}{$N$} & \multirow[c]{3}{*}{$R$} & \multirow[c]{3}{*}{$P$} & \multirow[c]{3}{*}{$C_p$} & \multicolumn{3}{c|}{\cite{valle2021}} & \multicolumn{3}{c|}{Strategy 1} & \multicolumn{7}{c|}{Strategy 2} & \multicolumn{7}{c|}{Strategy 3}\\
\cline{6-25}
&  &  &  &  & \multirow[c]{2}{*}{T(s)} & \multirow[c]{2}{*}{GAP} & \multirow[c]{2}{*}{UB} & \multicolumn{2}{c}{1st stage} & \multirow[c]{2}{*}{$|\Theta|$} & \multicolumn{2}{c}{1st stage} & \multicolumn{2}{c}{2nd stage} & \multirow[c]{2}{*}{$|\Theta|$} & \multirow[c]{2}{*}{$|S|$} & \multirow[c]{2}{*}{BLG} & \multicolumn{2}{c}{1st stage} & \multicolumn{2}{c}{2nd stage} & \multirow[c]{2}{*}{$|\Theta|$} & \multirow[c]{2}{*}{$|S|$} & \multirow[c]{2}{*}{BLG}\\
 &  & &  &  &  &  &  & \multicolumn{1}{c}{T(s)} & \multicolumn{1}{c}{GAP} &  & \multicolumn{1}{c}{T(s)} & \multicolumn{1}{c}{GAP} & \multicolumn{1}{c}{T(s)} & \multicolumn{1}{c}{GAP} &  &  &  & \multicolumn{1}{c}{T(s)} & \multicolumn{1}{c}{GAP} & \multicolumn{1}{c}{T(s)} & \multicolumn{1}{c}{GAP} &  &  & \\
\hline
\multirow{10}{*}{$50$} & \multirow{2}{*}{100} & \multirow{2}{*}{50} & 5 & 10 &  64.4 & -- & 7 &  78.4 & -- & 7 &   8.2 & -- &   0.1 & -- & 7 & 34 & 5 &   0.3 & -- &   0.0 & -- & 8 & 0 & 0\\
& &  & 10 & 5 &   1.3 & -- & 10 &   3.8 & -- & 10 &   2.0 & -- &   0.2 & -- & 10 & 34 & 0 &   0.3 & -- &   0.1 & -- & 10 & 0 & 0\\
\cline{2-25}
& \multirow{2}{*}{200} & \multirow{2}{*}{50} & 5 & 10 &   2.8 & -- & 12 &   2.7 & -- & 12 &   0.3 & -- &   0.1 & -- & 11 & 30 & 5 &   0.1 & -- &   0.1 & -- & 13 & 8 & 0\\
& &  & 10 & 5 &  95.8 & -- & 15 & 108.1 & -- & 15 &  33.9 & -- &   3.6 & -- & 15 & 30 & 7 &   0.6 & -- &   0.9 & -- & 13 & 8 & 4\\
\cline{2-25}
& \multirow{2}{*}{300} & \multirow{2}{*}{50} & 5 & 10 &   0.7 & -- & 16 &   1.1 & -- & 16 &   0.1 & -- &   0.1 & -- & 15 & 31 & 6 &   0.1 & -- &   0.1 & -- & 17 & 9 & 1\\
& &  & 10 & 5 &  45.4 & -- & 18 & 134.2 & -- & 18 &  40.2 & -- &   0.6 & -- & 17 & 31 & 9 &   0.2 & -- &   1.3 & -- & 17 & 9 & 2\\
\cline{2-25}
& \multirow{2}{*}{400} & \multirow{2}{*}{75} & 5 & 10 &   0.5 & -- & 19 &   0.3 & -- & 19 &   0.0 & -- &   0.0 & -- & 19 & 34 & 0 &   0.1 & -- &   0.0 & -- & 19 & 10 & 0\\
& &  & 10 & 5 &   2.5 & -- & 20 &   7.4 & -- & 20 &   3.5 & -- &   1.1 & -- & 19 & 34 & 5 &   0.1 & -- &   1.7 & -- & 19 & 10 & 2\\
\cline{2-25}
& \multirow{2}{*}{500} & \multirow{2}{*}{75} & 5 & 10 &   0.5 & -- & 20 &   0.1 & -- & 20 &   0.1 & -- &   0.1 & -- & 20 & 28 & 9 &   0.0 & -- &   0.1 & -- & 21 & 14 & 1\\
& &  & 10 & 5 &  19.3 & -- & 23 &  70.3 & -- & 23 &  17.5 & -- &   0.3 & -- & 21 & 28 & 7 &   0.1 & -- &   0.4 & -- & 21 & 14 & 3\\
\hline

\multirow{10}{*}{$100$} & \multirow{2}{*}{100} & \multirow{2}{*}{100} & 5 & 20 & TL & 31.16 & 9 & TL & 40.00 & 10 & TL & 32.33 &   0.2 & -- & 9 & 71 & 0 &   2.4 & -- &   0.1 & -- & 11 & 0 & 0\\
&  &  & 10 & 10 & TL &  9.09 & 11 & TL & 23.08 & 13 & TL &  9.02 &   3.6 & -- & 11 & 71 & 0 &  38.2 & -- &   0.2 & -- & 11 & 0 & 0\\
\cline{2-25}
& \multirow{2}{*}{200} & \multirow{2}{*}{100} & 5 & 20 & TL & 32.96 & 17 & TL & 30.16 & 17 & TL & 17.26 &   0.6 & -- & 15 & 59 & 16 &   0.3 & -- &   0.1 & -- & 20 & 8 & 0\\
&  &  & 10 & 10 & TL & 44.71 & 20 & TL & 47.00 & 21 & TL & 29.55 &   6.9 & -- & 18 & 59 & 3 &   0.8 & -- &   0.4 & -- & 20 & 8 & 0\\
\cline{2-25}
& \multirow{2}{*}{300} & \multirow{2}{*}{100} & 5 & 20 & TL & 16.11 & 20 & TL & 13.85 & 20 & 123.7 & -- &   0.2 & -- & 18 & 52 & 23 &   0.5 & -- &   0.1 & -- & 23 & 14 & 0\\
&  &  & 10 & 10 & TL & 31.77 & 24 & TL & 38.99 & 27 & TL & 26.47 &   6.4 & -- & 23 & 52 & 13 &   1.5 & -- &   1.3 & -- & 23 & 14 & 1\\
\cline{2-25}
& \multirow{2}{*}{400} & \multirow{2}{*}{100} & 5 & 20 &  98.3 & -- & 22 & 263.9 & -- & 22 &  41.7 & -- &   0.3 & -- & 22 & 60 & 27 &   0.5 & -- &   0.1 & -- & 25 & 16 & 0\\
& &  & 10 & 10 & TL & 30.87 & 29 & TL & 29.59 & 29 & TL & 11.82 &   3.6 & -- & 24 & 60 & 21 &   2.0 & -- &   3.1 & -- & 25 & 16 & 2\\
\cline{2-25}
& \multirow{2}{*}{500} & \multirow{2}{*}{100} & 5 & 20 &  75.3 & -- & 27 &  65.6 & -- & 27 &   7.8 & -- &   0.3 & -- & 26 & 61 & 20 &   0.9 & -- &   0.4 & -- & 29 & 23 & 0\\
& &  & 10 & 10 & TL & 21.50 & 33 & TL & 15.65 & 31 & TL &  5.75 &  13.1 & -- & 28 & 61 & 10 &   1.1 & -- &  16.0 & -- & 29 & 23 & 2\\
\hline

 \multirow{10}{*}{$150$} & \multirow{2}{*}{100} & \multirow{2}{*}{100} & 5 & 30 & TL & 37.22 & 13 & TL & 37.13 & 13 & TL & 35.60 &   0.7 & -- & 13 & 85 & 0 &  31.0 & -- &   0.2 & -- & 17 & 0 & 0\\
& &  & 10 & 15 & TL & 33.33 & 15 & TL & 37.50 & 16 & TL & 33.35 &   2.6 & -- & 15 & 85 & 0 & 255.4 & -- &   0.7 & -- & 17 & 0 & 0\\
\cline{2-25}
& \multirow{2}{*}{200} & \multirow{2}{*}{100} & 5 & 30 & TL & 40.09 & 21 & TL & 33.66 & 20 & TL & 24.06 &   0.7 & -- & 18 & 93 & 0 &   0.4 & -- &   0.1 & -- & 24 & 5 & 0\\
& &  & 10 & 15 & TL & 45.17 & 22 & TL & 51.44 & 25 & TL & 41.37 &  60.8 & -- & 22 & 93 & 4 &   1.3 & -- &   0.6 & -- & 24 & 5 & 0\\
\cline{2-25}
& \multirow{2}{*}{300} & \multirow{2}{*}{100} & 5 & 30 & TL & 26.46 & 27 & TL & 25.44 & 27 & TL & 12.90 &   0.9 & -- & 24 & 89 & 19 &   0.7 & -- &   1.3 & -- & 27 & 30 & 0\\
& &  & 10 & 15 & TL & 41.06 & 32 & TL & 45.65 & 35 & TL & 36.02 &   2.1 & -- & 31 & 89 & 0 &   4.8 & -- & TL & 12.28 & 27 & 30 & 5\\
\cline{2-25}
& \multirow{2}{*}{400} & \multirow{2}{*}{150} & 5 & 30 & TL & 24.43 & 28 & TL & 26.86 & 29 & TL & 11.16 &   0.4 & -- & 25 & 97 & 38 &   3.5 & -- &   0.2 & -- & 28 & 23 & 0\\
& &  & 10 & 15 & TL & 38.28 & 33 & TL & 42.13 & 36 & TL & 29.57 &  38.6 & -- & 30 & 97 & 27 & 186.7 & -- &  25.4 & -- & 28 & 23 & 2\\
\cline{2-25}
& \multirow{2}{*}{500} & \multirow{2}{*}{150} & 5 & 30 & TL & 23.29 & 35 & TL & 17.83 & 33 & TL &  2.67 &   0.4 & -- & 29 & 88 & 38 &   1.3 & -- &   0.8 & -- & 33 & 30 & 1\\
& &  & 10 & 15 & TL & 38.51 & 42 & TL & 36.02 & 41 & TL & 25.69 &  26.7 & -- & 36 & 88 & 25 &  12.7 & -- & 146.2 & -- & 33 & 30 & 5\\

\hline
\multirow{10}{*}{$200$} &\multirow{2}{*}{100} & \multirow{2}{*}{200} & 5 & 40 & TL & 47.85 & 14 & TL & 51.04 & 15 & TL & 42.29 &  16.1 & -- & 13 & 121 & 11 & TL & 11.53 &   0.4 & -- & 17 & 0 & 0\\
& &  & 10 & 20 & TL & 44.44 & 18 & TL & 41.18 & 17 & TL & 33.36 &   9.4 & -- & 15 & 121 & 0 & TL & 75.67 &   0.3 & -- & 60 & 0 & 0\\
\cline{2-25}
& \multirow{2}{*}{200} & \multirow{2}{*}{200} & 5 & 40 & TL & 49.08 & 23 & TL & 42.63 & 21 & TL & 38.83 &   1.5 & -- & 20 & 119 & 0 &   1.3 & -- &   0.3 & -- & 33 & 4 & 0\\
& &  & 10 & 20 & TL & 55.35 & 26 & TL & 76.45 & 51 & TL & 53.55 & TL &  9.24 & 26 & 119 & 11 &   4.3 & -- &   1.0 & -- & 33 & 4 & 0\\
\cline{2-25}
& \multirow{2}{*}{300} & \multirow{2}{*}{200} & 5 & 40 & TL & 41.67 & 28 & TL & 41.64 & 28 & TL & 34.55 &   1.2 & -- & 26 & 123 & 0 &   4.8 & -- &   0.3 & -- & 34 & 20 & 0\\
 & &  & 10 & 20 & TL & 80.21 & 82 & TL & 54.54 & 36 & TL & 42.83 & TL & 16.26 & 29 & 123 & 20 &  32.3 & -- &   3.4 & -- & 34 & 20 & 0\\
\cline{2-25}
& \multirow{2}{*}{400} & \multirow{2}{*}{200} & 5 & 40 & TL & 35.54 & 31 & TL & 39.76 & 34 & TL & 27.98 &   1.1 & -- & 29 & 128 & 0 &   9.3 & -- &   0.4 & -- & 34 & 27 & 0\\
 & &  & 10 & 20 & TL & 47.45 & 37 & TL & 50.72 & 41 & TL & 40.51 &   7.7 & -- & 34 & 128 & 0 & TL &  2.98 &  25.6 & -- & 34 & 27 & 0\\
\cline{2-25}
& \multirow{2}{*}{500} & \multirow{2}{*}{200} & 5 & 40 & TL & 30.70 & 40 & TL & 32.01 & 42 & TL & 23.70 &   3.0 & -- & 38 & 121 & 0 &   4.2 & -- &   0.7 & -- & 44 & 34 & 0\\
&  &  & 10 & 20 & TL & 39.24 & 45 & TL & 39.84 & 47 & TL & 34.24 & TL &  0.83 & 43 & 121 & 1 &  22.2 & -- &   7.3 & -- & 44 & 34 & 0\\

\hline
\multirow{10}{*}{$500$} & \multirow{2}{*}{100} & \multirow{2}{*}{200} & 5 & 100 & TL & 37.05 & 23 & TL & 33.70 & 22 & TL & 22.78 &   5.5 & -- & 19 & 303 & 0 & TL & 74.18 &   4.5 & -- & 80 & 0 & 0\\
& &  & 10 & 50 & TL & 63.94 & 40 & TL & 58.54 & 35 & TL & 38.99 &  31.5 & -- & 24 & 303 & 0 & TL &  &  &  & NF & 0 & \\
\cline{2-25}
& \multirow{2}{*}{200} & \multirow{2}{*}{200} & 5 & 100 & TL & 44.19 & 35 & TL & 43.37 & 36 & TL & 35.86 &  12.8 & -- & 32 & 285 & 0 &   3.1 & -- &   1.8 & -- & 63 & 16 & 0\\
&  &  & 10 & 50 & TL &  & NF & TL &  & NF & TL & 54.08 &  56.8 & -- & 44 & 285 & 0 &  15.3 & -- &  26.2 & -- & 63 & 16 & 0\\
\cline{2-25}
& \multirow{2}{*}{300} & \multirow{2}{*}{200} & 5 & 100 & TL & 74.94 & 93 & TL & 64.11 & 69 & TL & 39.66 &   9.2 & -- & 41 & 292 & 0 &   0.5 &  &  &  & NF & 41 & \\
&  &  & 10 & 50 & TL &  & NF & TL &  & NF & TL & 48.62 &  91.3 & -- & 48 & 292 & 0 &   1.0 &  &  &  & NF & 41 & \\
\cline{2-25}
& \multirow{2}{*}{400} & \multirow{2}{*}{200} & 5 & 100 & TL & 63.85 & 84 & TL & 55.66 & 71 & TL & 39.07 &  12.9 & -- & 52 & 275 & 0 &   3.6 & -- &   2.2 & -- & 68 & 75 & 0\\
&  &  & 10 & 50 & TL &  & NF & TL &  & NF & TL & 49.77 & 110.0 & -- & 62 & 275 & 0 &  56.6 & -- &  27.1 & -- & 68 & 75 & 0\\
\cline{2-25}
& \multirow{2}{*}{500} & \multirow{2}{*}{200} & 5 & 100 & TL & 52.04 & 77 & TL & 43.03 & 68 & TL & 31.99 &  11.6 & -- & 57 & 307 & 0 & 125.2 & -- &   2.2 & -- & 68 & 67 & 0\\
&  &  & 10 & 50 & TL &  & NF & TL &  & NF & TL & 35.68 & 262.4 & -- & 60 & 307 & 0 & TL &  1.77 &  36.2 & -- & 68 & 67 & 0\\

\hline
\multicolumn{3}{l}{\multirow{5}{*}{\textbf{Average:}}} & \multicolumn{2}{l}{$O = 50$} & \multicolumn{1}{r}{\textbf{ 23.3}} & \multicolumn{1}{r}{\textbf{  0.0}} & \multicolumn{1}{r}{\textbf{ 16.0}} & \multicolumn{1}{r}{\textbf{ 40.6}} & \multicolumn{1}{r}{\textbf{  0.0}} & \multicolumn{1}{r}{\textbf{ 16.0}} & \multicolumn{1}{r}{\textbf{ 10.6}} & \multicolumn{1}{r}{\textbf{  0.0}} & \multicolumn{1}{r}{\textbf{  0.6}} & \multicolumn{1}{r}{\textbf{  0.0}} & \multicolumn{1}{r}{\textbf{ 15.4}} & \multicolumn{1}{r}{\textbf{ 31.4}} & \multicolumn{1}{r}{\textbf{  5.3}} & \multicolumn{1}{r}{\textbf{  0.2}} & \multicolumn{1}{r}{\textbf{  0.0}} & \multicolumn{1}{r}{\textbf{  0.5}} & \multicolumn{1}{r}{\textbf{  0.0}} & \multicolumn{1}{r}{\textbf{ 15.8}} & \multicolumn{1}{r}{\textbf{  8.2}} & \multicolumn{1}{r}{\textbf{  1.3}}\\

\multicolumn{3}{l}{} & \multicolumn{2}{l}{$O = 100$} & \multicolumn{1}{r}{\textbf{   257.4}} & \multicolumn{1}{r}{\textbf{ 21.8}} & \multicolumn{1}{r}{\textbf{ 21.2}} & \multicolumn{1}{r}{\textbf{   272.9}} & \multicolumn{1}{r}{\textbf{ 23.8}} & \multicolumn{1}{r}{\textbf{ 21.7}} & \multicolumn{1}{r}{\textbf{   227.3}} & \multicolumn{1}{r}{\textbf{ 13.2}} & \multicolumn{1}{r}{\textbf{  3.5}} & \multicolumn{1}{r}{\textbf{  0.0}} & \multicolumn{1}{r}{\textbf{ 19.4}} & \multicolumn{1}{r}{\textbf{ 60.6}} & \multicolumn{1}{r}{\textbf{ 13.3}} & \multicolumn{1}{r}{\textbf{  4.8}} & \multicolumn{1}{r}{\textbf{  0.0}} & \multicolumn{1}{r}{\textbf{  2.2}} & \multicolumn{1}{r}{\textbf{  0.0}} & \multicolumn{1}{r}{\textbf{ 21.6}} & \multicolumn{1}{r}{\textbf{ 12.2}} & \multicolumn{1}{r}{\textbf{  0.5}}\\

\multicolumn{3}{l}{} & \multicolumn{2}{l}{$O = 150$} & \multicolumn{1}{r}{\textbf{   TL}} & \multicolumn{1}{r}{\textbf{ 34.8}} & \multicolumn{1}{r}{\textbf{ 26.8}} & \multicolumn{1}{r}{\textbf{   TL}} & \multicolumn{1}{r}{\textbf{ 35.4}} & \multicolumn{1}{r}{\textbf{ 27.5}} & \multicolumn{1}{r}{\textbf{   TL}} & \multicolumn{1}{r}{\textbf{ 25.2}} & \multicolumn{1}{r}{\textbf{ 13.4}} & \multicolumn{1}{r}{\textbf{  0.0}} & \multicolumn{1}{r}{\textbf{ 24.3}} & \multicolumn{1}{r}{\textbf{ 90.4}} & \multicolumn{1}{r}{\textbf{ 15.1}} & \multicolumn{1}{r}{\textbf{ 49.8}} & \multicolumn{1}{r}{\textbf{  0.0}} & \multicolumn{1}{r}{\textbf{ 47.6}} & \multicolumn{1}{r}{\textbf{  1.2}} & \multicolumn{1}{r}{\textbf{ 25.8}} & \multicolumn{1}{r}{\textbf{ 17.6}} & \multicolumn{1}{r}{\textbf{  1.3}}\\

\multicolumn{3}{l}{} & \multicolumn{2}{l}{$O = 200$} & \multicolumn{1}{r}{\textbf{   TL}} & \multicolumn{1}{r}{\textbf{ 47.2}} & \multicolumn{1}{r}{\textbf{ 34.4}} & \multicolumn{1}{r}{\textbf{   TL}} & \multicolumn{1}{r}{\textbf{ 47.0}} & \multicolumn{1}{r}{\textbf{ 33.2}} & \multicolumn{1}{r}{\textbf{   TL}} & \multicolumn{1}{r}{\textbf{ 37.2}} & \multicolumn{1}{r}{\textbf{ 94.0}} & \multicolumn{1}{r}{\textbf{  2.6}} & \multicolumn{1}{r}{\textbf{ 27.3}} & \multicolumn{1}{r}{\textbf{   122.4}} & \multicolumn{1}{r}{\textbf{  4.3}} & \multicolumn{1}{r}{\textbf{ 97.8}} & \multicolumn{1}{r}{\textbf{  9.0}} & \multicolumn{1}{r}{\textbf{  4.0}} & \multicolumn{1}{r}{\textbf{  0.0}} & \multicolumn{1}{r}{\textbf{ 36.7}} & \multicolumn{1}{r}{\textbf{ 17.0}} & \multicolumn{1}{r}{\textbf{  0.0}}\\

\multicolumn{3}{l}{} & \multicolumn{2}{l}{$O = 500$} &  \multicolumn{1}{r}{\textbf{   TL}} & \multicolumn{1}{r}{\textbf{ 56.0}} & \multicolumn{1}{r}{\textbf{ 58.7}} & \multicolumn{1}{r}{\textbf{   TL}} & \multicolumn{1}{r}{\textbf{ 49.7}} & \multicolumn{1}{r}{\textbf{ 50.2}} & \multicolumn{1}{r}{\textbf{   TL}} & \multicolumn{1}{r}{\textbf{ 39.6}} & \multicolumn{1}{r}{\textbf{ 60.4}} & \multicolumn{1}{r}{\textbf{  0.0}} & \multicolumn{1}{r}{\textbf{ 43.9}} & \multicolumn{1}{r}{\textbf{   292.4}} & \multicolumn{1}{r}{\textbf{  0.0}} & \multicolumn{1}{r}{\textbf{   110.5}} & \multicolumn{1}{r}{\textbf{ 10.9}} & \multicolumn{1}{r}{\textbf{ 14.3}} & \multicolumn{1}{r}{\textbf{  0.0}} & \multicolumn{1}{r}{\textbf{ 68.3}} & \multicolumn{1}{r}{\textbf{ 39.8}} & \multicolumn{1}{r}{\textbf{  0.0}}

\end{tabular}
}
\caption{$O = 50, 100, 150, 200, 500$, no rack reduction}
\label{nored}
\end{table}


\begin{table}[!htpb]
\centering
{\tiny
\renewcommand{\tabcolsep}{0.8mm} 
\renewcommand{\arraystretch}{1.8} 
\begin{tabular}{|c|c|c|c|c|rrr|rrr|rrrrrrr|rrrrrrr|}
\hline
 \multirow[c]{3}{*}{$O$} &  \multirow[c]{3}{*}{$N$} & \multirow[c]{3}{*}{$R$} & \multirow[c]{3}{*}{$P$} & \multirow[c]{3}{*}{$C_p$} & \multicolumn{3}{c|}{\cite{valle2021}} & \multicolumn{3}{c|}{Strategy 1} & \multicolumn{7}{c|}{Strategy 2} & \multicolumn{7}{c|}{Strategy 3}\\
\cline{6-25}
&  &  &  &  & \multirow[c]{2}{*}{T(s)} & \multirow[c]{2}{*}{GAP} & \multirow[c]{2}{*}{UB} & \multicolumn{2}{c}{1st stage} & \multirow[c]{2}{*}{$|\Theta|$} & \multicolumn{2}{c}{1st stage} & \multicolumn{2}{c}{2nd stage} & \multirow[c]{2}{*}{$|\Theta|$} & \multirow[c]{2}{*}{$|S|$} & \multirow[c]{2}{*}{BLG} & \multicolumn{2}{c}{1st stage} & \multicolumn{2}{c}{2nd stage} & \multirow[c]{2}{*}{$|\Theta|$} & \multirow[c]{2}{*}{$|S|$} & \multirow[c]{2}{*}{BLG}\\
 &  & &  &  &  &  &  & \multicolumn{1}{c}{T(s)} & \multicolumn{1}{c}{GAP} &  & \multicolumn{1}{c}{T(s)} & \multicolumn{1}{c}{GAP} & \multicolumn{1}{c}{T(s)} & \multicolumn{1}{c}{GAP} &  &  &  & \multicolumn{1}{c}{T(s)} & \multicolumn{1}{c}{GAP} & \multicolumn{1}{c}{T(s)} & \multicolumn{1}{c}{GAP} &  &  & \\
\hline

\multirow{10}{*}{$50$} & \multirow{2}{*}{100} & \multirow{2}{*}{50} & 5 & 10 &   0.1 & -- & 8 &   0.0 & -- & 8 &   0.0 & -- &   0.1 & -- & 8 & 34 & 0 &   0.1 & -- &   0.0 & -- & 9 & 0 & 0\\
&  &  & 10 & 5 &   0.4 & -- & 11 &   0.5 & -- & 11 &   0.1 & -- &   0.3 & -- & 11 & 34 & 0 &   0.3 & -- &   0.1 & -- & 10 & 0 & 0\\
\cline{2-25}
& \multirow{2}{*}{200} & \multirow{2}{*}{50} & 5 & 10 &   0.3 & -- & 13 &   0.2 & -- & 13 &   0.0 & -- &   0.1 & -- & 11 & 30 & 6 &   0.0 & -- &   0.0 & -- & 13 & 8 & 0\\
&  &  & 10 & 5 &   7.3 & -- & 16 &  23.4 & -- & 16 &   1.5 & -- &   5.0 & -- & 16 & 30 & 3 &   0.1 & -- &   1.5 & -- & 14 & 8 & 2\\
\cline{2-25}
& \multirow{2}{*}{300} & \multirow{2}{*}{50} & 5 & 10 &   0.2 & -- & 19 &   0.1 & -- & 19 &   0.0 & -- &   0.1 & -- & 19 & 31 & 0 &   0.0 & -- &   0.0 & -- & 19 & 9 & 0\\
&  &  & 10 & 5 &   1.3 & -- & 20 &   3.0 & -- & 20 &   0.1 & -- &   1.0 & -- & 19 & 31 & 2 &   0.5 & -- &   1.4 & -- & 19 & 9 & 1\\
\cline{2-25}
& \multirow{2}{*}{400} & \multirow{2}{*}{75} & 5 & 10 &   0.3 & -- & 20 &   0.2 & -- & 20 &   0.0 & -- &   0.1 & -- & 20 & 34 & 1 &   0.0 & -- &   0.1 & -- & 20 & 10 & 1\\
&  &  & 10 & 5 &   0.3 & -- & 23 &   0.5 & -- & 23 &   0.1 & -- &   0.4 & -- & 21 & 34 & 3 &   0.1 & -- &   0.1 & -- & 23 & 10 & 0\\
\cline{2-25}
& \multirow{2}{*}{500} & \multirow{2}{*}{75} & 5 & 10 &   0.2 & -- & 22 &   0.1 & -- & 22 &   0.0 & -- &   0.1 & -- & 22 & 28 & 0 &   0.0 & -- &   0.1 & -- & 22 & 14 & 0\\
&  &  & 10 & 5 &   1.0 & -- & 25 &   1.9 & -- & 25 &   0.6 & -- &   0.5 & -- & 22 & 28 & 7 &   0.1 & -- &   0.6 & -- & 22 & 14 & 2\\

\hline
\multirow{10}{*}{$100$} & \multirow{2}{*}{100} & \multirow{2}{*}{100} & 5 & 20 &   0.4 & -- & 9 &   0.2 & -- & 9 &   0.0 & -- &   0.3 & -- & 9 & 71 & 0 &   0.4 & -- &   0.1 & -- & 13 & 0 & 0\\
&  &  & 10 & 10 &   5.8 & -- & 12 &   8.9 & -- & 12 &   0.8 & -- &   9.0 & -- & 12 & 71 & 0 &   1.6 & -- &   0.3 & -- & 13 & 0 & 0\\
\cline{2-25}
& \multirow{2}{*}{200} & \multirow{2}{*}{100} & 5 & 20 &   1.6 & -- & 16 &   0.7 & -- & 16 &   0.2 & -- &   0.3 & -- & 16 & 59 & 0 &   0.2 & -- &   0.1 & -- & 23 & 8 & 0\\
&  &  & 10 & 10 & TL & 19.79 & 20 & TL & 13.71 & 19 & TL &  3.80 &  34.6 & -- & 18 & 59 & 3 &   0.5 & -- &   0.5 & -- & 24 & 8 & 0\\
\cline{2-25}
& \multirow{2}{*}{300} & \multirow{2}{*}{100} & 5 & 20 &   3.7 & -- & 22 &   1.9 & -- & 22 &   0.1 & -- &   0.3 & -- & 22 & 52 & 0 &   0.1 & -- &   0.1 & -- & 26 & 14 & 0\\
&  &  & 10 & 10 & 206.4 & -- & 24 & 241.3 & -- & 24 & 259.7 & -- &  12.8 & -- & 23 & 52 & 5 &   0.4 & -- &   5.2 & -- & 26 & 14 & 1\\
\cline{2-25}
& \multirow{2}{*}{400} & \multirow{2}{*}{100} & 5 & 20 &   0.7 & -- & 25 &   0.2 & -- & 25 &   0.1 & -- &   0.2 & -- & 25 & 60 & 0 &   0.1 & -- &   0.1 & -- & 28 & 16 & 0\\
&  &  & 10 & 10 & 225.4 & -- & 27 & 154.4 & -- & 27 &   0.2 & -- &   5.1 & -- & 26 & 60 & 8 &   0.3 & -- &   0.8 & -- & 29 & 16 & 0\\
\cline{2-25}
& \multirow{2}{*}{500} & \multirow{2}{*}{100} & 5 & 20 &   0.4 & -- & 29 &   0.3 & -- & 29 &   0.1 & -- &   0.3 & -- & 29 & 61 & 0 &   0.1 & -- &   0.1 & -- & 31 & 23 & 0\\
&  &  & 10 & 10 & TL &  8.99 & 32 & TL &  6.25 & 32 &   2.3 & -- &  11.8 & -- & 29 & 61 & 11 &   0.2 & -- &   5.1 & -- & 31 & 23 & 1\\

\hline
\multirow{10}{*}{$150$} & \multirow{2}{*}{100} & \multirow{2}{*}{100} & 5 & 30 &  81.7 & -- & 13 &   8.1 & -- & 13 &   0.2 & -- &   6.2 & -- & 12 & 85 & 3 &   1.0 & -- &   0.1 & -- & 21 & 0 & 0\\
&  &  & 10 & 15 & 115.4 & -- & 14 & 221.5 & -- & 14 &  78.5 & -- &  45.1 & -- & 13 & 85 & 4 &   3.1 & -- &   1.6 & -- & 19 & 0 & 0\\
\cline{2-25}
& \multirow{2}{*}{200} & \multirow{2}{*}{100} & 5 & 30 &  28.1 & -- & 20 &   1.4 & -- & 20 &   0.1 & -- &   0.5 & -- & 20 & 93 & 0 &   0.4 & -- &   0.2 & -- & 27 & 5 & 0\\
&  &  & 10 & 15 & TL & 13.64 & 22 & TL & 13.64 & 22 & 108.3 & -- &  79.3 & -- & 19 & 93 & 17 &   1.1 & -- &   0.6 & -- & 30 & 5 & 0\\
\cline{2-25}
& \multirow{2}{*}{300} & \multirow{2}{*}{100} & 5 & 30 &  30.6 & -- & 26 &   7.4 & -- & 26 &   0.4 & -- &   0.6 & -- & 26 & 89 & 0 &   0.2 & -- &   0.7 & -- & 29 & 30 & 0\\
&  &  & 10 & 15 & TL & 28.71 & 36 & TL & 16.13 & 31 & TL &  3.70 & TL & 22.47 & 27 & 89 & 20 &   0.4 & -- &   3.8 & -- & 32 & 30 & 0\\
\cline{2-25}
& \multirow{2}{*}{400} & \multirow{2}{*}{150} & 5 & 30 &   5.1 & -- & 29 &   1.0 & -- & 29 &   0.1 & -- &   0.3 & -- & 29 & 97 & 0 &   0.3 & -- &   0.2 & -- & 33 & 23 & 0\\
&  &  & 10 & 15 & TL &  9.68 & 31 & TL &  9.68 & 31 & 261.4 & -- &  21.2 & -- & 28 & 97 & 17 &   0.8 & -- &  14.0 & -- & 32 & 23 & 0\\
\cline{2-25}
& \multirow{2}{*}{500} & \multirow{2}{*}{150} & 5 & 30 &   5.3 & -- & 34 &   1.5 & -- & 34 &   0.1 & -- &   0.7 & -- & 34 & 88 & 0 &   0.2 & -- &   0.4 & -- & 38 & 30 & 0\\
&  &  & 10 & 15 & TL & 12.19 & 39 & TL & 96.73 & 43 & TL &  2.84 &  25.8 & -- & 35 & 88 & 22 &   0.5 & -- & 108.6 & -- & 37 & 30 & 1\\

\hline
\multirow{10}{*}{$200$} & \multirow{2}{*}{100} & \multirow{2}{*}{200} & 5 & 40 &  84.3 & -- & 12 &   4.7 & -- & 12 &   0.3 & -- &   1.4 & -- & 12 & 121 & 0 &   2.8 & -- &   0.4 & -- & 22 & 0 & 0\\
&  &  & 10 & 20 & TL & 15.08 & 15 & TL & 18.19 & 16 & TL & 18.74 &  20.2 & -- & 16 & 121 & 0 &   8.5 & -- &   1.4 & -- & 25 & 0 & 0\\
\cline{2-25}
& \multirow{2}{*}{200} & \multirow{2}{*}{200} & 5 & 40 & 109.4 & -- & 20 & 103.3 & -- & 20 &   0.5 & -- &   3.5 & -- & 20 & 119 & 0 &   1.2 & -- &   0.3 & -- & 36 & 4 & 0\\
&  &  & 10 & 20 & TL & 34.48 & 29 & TL & 29.52 & 27 & TL & 20.20 & TL &  0.84 & 24 & 119 & 1 &   2.4 & -- &   1.2 & -- & 36 & 4 & 0\\
\cline{2-25}
& \multirow{2}{*}{300} & \multirow{2}{*}{200} & 5 & 40 &  48.6 & -- & 26 &   8.7 & -- & 26 &   0.2 & -- &   1.7 & -- & 26 & 123 & 0 &   0.8 & -- &   0.4 & -- & 39 & 20 & 0\\
&  &  & 10 & 20 & TL & 20.59 & 34 & TL & 18.18 & 33 & TL &  3.57 & TL & 17.07 & 28 & 123 & 21 &   1.4 & -- &   3.2 & -- & 39 & 20 & 0\\
\cline{2-25}
& \multirow{2}{*}{400} & \multirow{2}{*}{200} & 5 & 40 &   4.0 & -- & 31 &   2.9 & -- & 31 &   0.2 & -- &   1.2 & -- & 31 & 128 & 0 &   0.7 & -- &   0.6 & -- & 39 & 27 & 0\\
&  &  & 10 & 20 & TL & 15.22 & 38 & TL &  5.71 & 35 &   8.0 & -- & TL &  4.69 & 33 & 128 & 6 &   1.8 & -- &   4.8 & -- & 37 & 27 & 0\\
\cline{2-25}
& \multirow{2}{*}{500} & \multirow{2}{*}{200} & 5 & 40 &  12.1 & -- & 42 &   3.4 & -- & 42 &   0.2 & -- &   0.6 & -- & 42 & 121 & 0 &   0.5 & -- &   0.4 & -- & 52 & 34 & 0\\
&  &  & 10 & 20 & TL & 17.67 & 49 & TL &  8.89 & 45 &  96.8 & -- & TL & 14.66 & 41 & 121 & 20 &   1.0 & -- &   4.1 & -- & 50 & 34 & 0\\

\hline
\multirow{10}{*}{$500$} & \multirow{2}{*}{100} & \multirow{2}{*}{200} & 5 & 100 & 187.3 & -- & 21 &   5.9 & -- & 21 &   0.3 & -- &   6.0 & -- & 21 & 303 & 0 &  35.3 & -- &   3.8 & -- & 32 & 0 & 0\\
&  &  & 10 & 50 & TL &  & NF & TL &  9.44 & 22 & 270.6 & -- &  28.7 & -- & 21 & 303 & 0 & TL &  &  &  & NF & 0 & \\
\cline{2-25}
& \multirow{2}{*}{200} & \multirow{2}{*}{200} & 5 & 100 & TL & 22.56 & 36 & 247.3 & -- & 29 &   8.8 & -- & 127.2 & -- & 29 & 285 & 0 &   3.6 & -- &   1.9 & -- & 70 & 16 & 0\\
&  &  & 10 & 50 & TL &  & NF & TL & 32.16 & 40 & TL & 27.30 &  79.4 & -- & 38 & 285 & 0 &  14.1 & -- &   9.2 & -- & 71 & 16 & 0\\
\cline{2-25}
& \multirow{2}{*}{300} & \multirow{2}{*}{200} & 5 & 100 & TL & 16.86 & 44 &  29.6 & -- & 38 &   1.4 & -- &  10.8 & -- & 38 & 292 & 0 &   0.5 &  &  &  & NF & 41 & \\
&  &  & 10 & 50 & TL &  & NF & TL & 34.43 & 54 & TL & 15.40 & TL & 18.49 & 42 & 292 & 54 &   0.9 &  &  &  & NF & 41 & \\
\cline{2-25}
& \multirow{2}{*}{400} & \multirow{2}{*}{200} & 5 & 100 & TL & 10.07 & 56 &  45.4 & -- & 51 &   2.4 & -- &  17.8 & -- & 51 & 275 & 0 &   1.7 & -- &   2.9 & -- & 70 & 75 & 0\\
&  &  & 10 & 50 & TL &  & NF & TL &  & NF & TL & 18.05 & 133.6 & -- & 56 & 275 & 0 &  16.8 & -- &  38.9 & -- & 73 & 75 & 0\\
\cline{2-25}
& \multirow{2}{*}{500} & \multirow{2}{*}{200} & 5 & 100 & TL &  9.96 & 58 & 244.5 & -- & 53 &   2.4 & -- &  12.4 & -- & 52 & 307 & 0 &   2.8 & -- &   2.8 & -- & 72 & 67 & 0\\
& &  & 10 & 50 & TL &  & NF & TL & 18.46 & 65 & TL &  5.36 & TL & 10.75 & 56 & 307 & 33 &   4.1 & -- &  38.0 & -- & 73 & 67 & 0\\

\hline

 \multicolumn{3}{l}{\multirow{5}{*}{\textbf{Average:}}} & \multicolumn{2}{l}{$O = 50$} & \multicolumn{1}{r}{\textbf{  1.1}} & \multicolumn{1}{r}{\textbf{  0.0}} & \multicolumn{1}{r}{\textbf{ 17.7}} & \multicolumn{1}{r}{\textbf{  3.0}} & \multicolumn{1}{r}{\textbf{  0.0}} & \multicolumn{1}{r}{\textbf{ 17.7}} & \multicolumn{1}{r}{\textbf{  0.3}} & \multicolumn{1}{r}{\textbf{  0.0}} & \multicolumn{1}{r}{\textbf{  0.8}} & \multicolumn{1}{r}{\textbf{  0.0}} & \multicolumn{1}{r}{\textbf{ 16.9}} & \multicolumn{1}{r}{\textbf{ 31.4}} & \multicolumn{1}{r}{\textbf{  2.2}} & \multicolumn{1}{r}{\textbf{  0.1}} & \multicolumn{1}{r}{\textbf{  0.0}} & \multicolumn{1}{r}{\textbf{  0.4}} & \multicolumn{1}{r}{\textbf{  0.0}} & \multicolumn{1}{r}{\textbf{ 17.1}} & \multicolumn{1}{r}{\textbf{  8.2}} & \multicolumn{1}{r}{\textbf{  0.6}}\\

\multicolumn{3}{l}{} & \multicolumn{2}{l}{$O = 100$} &  \multicolumn{1}{r}{\textbf{   104.4}} & \multicolumn{1}{r}{\textbf{  2.9}} & \multicolumn{1}{r}{\textbf{ 21.6}} & \multicolumn{1}{r}{\textbf{   100.8}} & \multicolumn{1}{r}{\textbf{  2.0}} & \multicolumn{1}{r}{\textbf{ 21.5}} & \multicolumn{1}{r}{\textbf{ 56.3}} & \multicolumn{1}{r}{\textbf{  0.4}} & \multicolumn{1}{r}{\textbf{  7.5}} & \multicolumn{1}{r}{\textbf{  0.0}} & \multicolumn{1}{r}{\textbf{ 20.9}} & \multicolumn{1}{r}{\textbf{ 60.6}} & \multicolumn{1}{r}{\textbf{  2.7}} & \multicolumn{1}{r}{\textbf{  0.4}} & \multicolumn{1}{r}{\textbf{  0.0}} & \multicolumn{1}{r}{\textbf{  1.2}} & \multicolumn{1}{r}{\textbf{  0.0}} & \multicolumn{1}{r}{\textbf{ 24.4}} & \multicolumn{1}{r}{\textbf{ 12.2}} & \multicolumn{1}{r}{\textbf{  0.2}}\\

\multicolumn{3}{l}{} & \multicolumn{2}{l}{$O = 150$} &  \multicolumn{1}{r}{\textbf{   146.6}} & \multicolumn{1}{r}{\textbf{  6.4}} & \multicolumn{1}{r}{\textbf{ 26.4}} & \multicolumn{1}{r}{\textbf{   144.1}} & \multicolumn{1}{r}{\textbf{ 13.6}} & \multicolumn{1}{r}{\textbf{ 26.3}} & \multicolumn{1}{r}{\textbf{   104.9}} & \multicolumn{1}{r}{\textbf{  0.7}} & \multicolumn{1}{r}{\textbf{ 48.0}} & \multicolumn{1}{r}{\textbf{  2.2}} & \multicolumn{1}{r}{\textbf{ 24.3}} & \multicolumn{1}{r}{\textbf{ 90.4}} & \multicolumn{1}{r}{\textbf{  8.3}} & \multicolumn{1}{r}{\textbf{  0.8}} & \multicolumn{1}{r}{\textbf{  0.0}} & \multicolumn{1}{r}{\textbf{ 13.0}} & \multicolumn{1}{r}{\textbf{  0.0}} & \multicolumn{1}{r}{\textbf{ 29.8}} & \multicolumn{1}{r}{\textbf{ 17.6}} & \multicolumn{1}{r}{\textbf{  0.1}}\\

\multicolumn{3}{l}{} & \multicolumn{2}{l}{$O = 200$} &  \multicolumn{1}{r}{\textbf{   175.9}} & \multicolumn{1}{r}{\textbf{ 10.3}} & \multicolumn{1}{r}{\textbf{ 29.6}} & \multicolumn{1}{r}{\textbf{   162.3}} & \multicolumn{1}{r}{\textbf{  8.0}} & \multicolumn{1}{r}{\textbf{ 28.7}} & \multicolumn{1}{r}{\textbf{   100.6}} & \multicolumn{1}{r}{\textbf{  4.3}} & \multicolumn{1}{r}{\textbf{   122.8}} & \multicolumn{1}{r}{\textbf{  3.7}} & \multicolumn{1}{r}{\textbf{ 27.3}} & \multicolumn{1}{r}{\textbf{   122.4}} & \multicolumn{1}{r}{\textbf{  4.8}} & \multicolumn{1}{r}{\textbf{  2.1}} & \multicolumn{1}{r}{\textbf{  0.0}} & \multicolumn{1}{r}{\textbf{  1.7}} & \multicolumn{1}{r}{\textbf{  0.0}} & \multicolumn{1}{r}{\textbf{ 37.5}} & \multicolumn{1}{r}{\textbf{ 17.0}} & \multicolumn{1}{r}{\textbf{  0.0}}\\

\multicolumn{3}{l}{} & \multicolumn{2}{l}{$O = 500$} &  \multicolumn{1}{r}{\textbf{   288.7}} & \multicolumn{1}{r}{\textbf{ 11.9}} & \multicolumn{1}{r}{\textbf{ 43.0}} & \multicolumn{1}{r}{\textbf{   207.3}} & \multicolumn{1}{r}{\textbf{ 10.5}} & \multicolumn{1}{r}{\textbf{ 41.4}} & \multicolumn{1}{r}{\textbf{   148.6}} & \multicolumn{1}{r}{\textbf{  6.6}} & \multicolumn{1}{r}{\textbf{   101.6}} & \multicolumn{1}{r}{\textbf{  2.9}} & \multicolumn{1}{r}{\textbf{ 40.4}} & \multicolumn{1}{r}{\textbf{   292.4}} & \multicolumn{1}{r}{\textbf{  8.7}} & \multicolumn{1}{r}{\textbf{ 38.0}} & \multicolumn{1}{r}{\textbf{  0.0}} & \multicolumn{1}{r}{\textbf{ 13.9}} & \multicolumn{1}{r}{\textbf{  0.0}} & \multicolumn{1}{r}{\textbf{ 65.9}} & \multicolumn{1}{r}{\textbf{ 39.8}} & \multicolumn{1}{r}{\textbf{  0.0}}

\end{tabular}
}
\caption{$O = 50, 100, 150, 200, 500$, rack reduction}
\label{red}
\end{table}


\subsection{Comparison results}

Table~\ref{red} give detailed  results for the same problems as considered in Table~\ref{nored}, but now with rack reduction. Recall from above that our rack reduction procedure reduces the number of racks that need to be considered and so we would expect a decrease in computation time as a result of this reduction.

Comparing the results for~\cite{valle2021} and Strategy 1 we can conclude that with rack  reduction as problem size increases Strategy 1 performs better than the approach 
of~\cite{valle2021}. Here for $O \geq 100$ the results from Strategy 1 are better (on average) than those 
from~\cite{valle2021}.
For Strategies 2 and 3 we do see a decrease in computation time (considering both stages together). when we have rack reduction.

For comparison purposes 
Table~\ref{order_compare} shows in the first two groups of results  the summary rows for each value of $O$ as taken from 
Table~\ref{nored} and Table~\ref{red}. The final two groups of results in Table~\ref{order_compare} are summary results for the combined approach as discussed in Section~\ref{jebcomb} above. Detailed results for this combined approach are not given here for space reasons. 
The average rows in Table~\ref{order_compare}  gives the average value over all values of $O$. These are the averages over five different values for $O$, with ten cases being considered for each value of $O$.
Considering Table~\ref{order_compare} it seems reasonable to conclude that:
\begin{compactitem}
\item For both~\cite{valle2021} and Strategy 1, where we fulfil all orders at the first-stage, rack reduction has a very positive effect resulting (on average) in finding solutions requiring less racks quicker. 
\item For Strategy 2 and Strategy 3 where we fulfil some orders at the first-stage, some at the second-stage, we do not see such a marked improvement in performance in terms of the number of racks and the backlog as a result of rack reduction. However computation time (taking both stages into account) is (approximately) reduced by a factor of two for Strategy 2 and by a factor of five for Strategy 3.
\end{compactitem}


\noindent Recall from the discussion in Section~\ref{jebcomb} that the assumption behind  our two stage approach was that it will be computationally more effective to solve two smaller problems as compared with solving one larger (combined) problem.
Considering the combined approach results shown in Table~\ref{order_compare} for the two strategies (Strategies 2 and 3) that use two stages this assumption is justified. The average values (averaged over 50 instances) show that:
\begin{compactitem}
 \item the combined approach uses more racks and has a higher backlog; both when we have no rack reduction and when we have rack reduction
\item in all cases except for no rack reduction and Strategy 2 the combined approach requires more computation time as compared to the total computation time for stage 1 and stage 2. 
\end{compactitem}


\begin{table}[!htpb]

\centering
{\tiny
\renewcommand{\tabcolsep}{1mm} 
\renewcommand{\arraystretch}{1.8} 
\begin{tabular}{|r|c|rrr|rrr|rrrrrrr|rrrrrrr|}
\hline
\multirow[c]{3}{*}{Case} &
\multirow[c]{3}{*}{$O$} & \multicolumn{3}{c|}{\cite{valle2021}} & \multicolumn{3}{c|}{Strategy 1} & \multicolumn{7}{c|}{Strategy 2} & \multicolumn{7}{c|}{Strategy 3}\\
\cline{3-22}
 & &  \multirow[c]{2}{*}{T(s)} & \multirow[c]{2}{*}{GAP} & \multirow[c]{2}{*}{UB} & \multicolumn{2}{c}{1st stage} & \multirow[c]{2}{*}{$|\Theta|$} & \multicolumn{2}{c}{1st stage} & \multicolumn{2}{c}{2nd stage} & \multirow[c]{2}{*}{$|\Theta|$} & \multirow[c]{2}{*}{$|S|$} & \multirow[c]{2}{*}{BLG} & \multicolumn{2}{c}{1st stage} & \multicolumn{2}{c}{2nd stage} & \multirow[c]{2}{*}{$|\Theta|$} & \multirow[c]{2}{*}{$|S|$} & \multirow[c]{2}{*}{BLG}\\
 & &  &  &  & \multicolumn{1}{c}{T(s)} & \multicolumn{1}{c}{GAP} &  & \multicolumn{1}{c}{T(s)} & \multicolumn{1}{c}{GAP} & \multicolumn{1}{c}{T(s)} & \multicolumn{1}{c}{GAP} &  &  &  & \multicolumn{1}{c}{T(s)} & \multicolumn{1}{c}{GAP} & \multicolumn{1}{c}{T(s)} & \multicolumn{1}{c}{GAP} &  &  & \\
\hline

No & 50 & 23.3   &  0.0   & 16.0   & 40.6   &  0.0   & 16.0   & 10.6   &  0.0   &  0.6   &  0.0   & 15.4   & 31.4   &  5.3   &  0.2   &  0.0   &  0.5   &  0.0   & 15.8   &  8.2   &  1.3  \\

rack & 100 &   257.4   & 21.8   & 21.2   &   272.9   & 23.8   & 21.7   &   227.3   & 13.2   &  3.5   &  0.0   & 19.4   & 60.6   & 13.3   &  4.8   &  0.0   &  2.2   &  0.0   & 21.6   & 12.2   &  0.5  \\

reduction & 150 &   TL   & 34.8   & 26.8   &   TL   & 35.4   & 27.5   &   TL   & 25.2   & 13.4   &  0.0   & 24.3   & 90.4   & 15.1   & 49.8   &  0.0   & 47.6   &  1.2   & 25.8   & 17.6   &  1.3  \\

& 200 &   TL   & 47.2   & 34.4   &   TL   & 47.0   & 33.2   &   TL   & 37.2   & 94.0   &  2.6   & 27.3   &   122.4   &  4.3   & 97.8   &  9.0   &  4.0   &  0.0   & 36.7   & 17.0   &  0.0  \\

& 500 &   TL   & 56.0   & 58.7   &   TL   & 49.7   & 50.2   &   TL   & 39.6   & 60.4   &  0.0   & 43.9   &   292.4   &  0.0   &   110.5   & 10.9   & 14.3   &  0.0   & 68.3   & 39.8   &  0.0  \\

 & \textbf{ Average }	&	\textbf{236.1}	&	\textbf{32.0}	&	\textbf{31.4	}&	\textbf{242.7}	&	\textbf{31.2	}&	\textbf{29.7}	&	\textbf{227.6}	&	\textbf{23.0}	&	\textbf{34.4}	&	\textbf{0.5}	&	\textbf{26.1}	&	\textbf{119.4}	&	\textbf{7.6}	&	\textbf{52.6}	&	\textbf{4.0}&	\textbf{13.7}	&	\textbf{0.2}	&	\textbf{33.6}	&	\textbf{19.0}	&	\textbf{0.6}	\\

\hline

Rack & 50 &   1.1   &   0.0   &  17.7   &   3.0   &   0.0   &  17.7   &   0.3   &   0.0   &   0.8   &   0.0   &  16.9   &  31.4   &   2.2   &   0.1   &   0.0   &   0.4   &   0.0   &  17.1   &   8.2   &   0.6  \\

reduction & 100 &  104.4   &   2.9   &  21.6   &  100.8   &   2.0   &  21.5   &  56.3   &   0.4   &   7.5   &   0.0   &  20.9   &  60.6   &   2.7   &   0.4   &   0.0   &   1.2   &   0.0   &  24.4   &  12.2   &   0.2  \\

& 150 &  146.6   &   6.4   &  26.4   &  144.1   &  13.6   &  26.3   &  104.9   &   0.7   &  48.0   &   2.2   &  24.3   &  90.4   &   8.3   &   0.8   &   0.0   &  13.0   &   0.0   &  29.8   &  17.6   &   0.1  \\

& 200 &  175.9   &  10.3   &  29.6   &  162.3   &   8.0   &  28.7   &  100.6   &   4.3   &  122.8   &   3.7   &  27.3   &  122.4   &   4.8   &   2.1   &   0.0   &   1.7   &   0.0   &  37.5   &  17.0   &   0.0  \\

& 500 &  288.7   &  11.9   &  43.0   &  207.3   &  10.5   &  41.4   &  148.6   &   6.6   &  101.6   &   2.9   &  40.4   &  292.4   &   8.7   &  38.0   &   0.0   &  13.9   &   0.0   &  65.9   &  39.8   &   0.0  \\
& \textbf{Average} 	&	\textbf{143.3}	&	\textbf{6.3}	&	\textbf{27.7}	&	\textbf{123.5}	&	\textbf{6.8}&	\textbf{27.1}	&	\textbf{82.1}	&	\textbf{2.4}	&	\textbf{56.1}	&	\textbf{1.8}&	\textbf{26.0}	&	\textbf{119.4}	&	\textbf{5.3}	&	\textbf{8.3}	&	\textbf{0.0}&	\textbf{6.0}	&	\textbf{0.0}	&	\textbf{34.9}&	\textbf{19.0}	&	\textbf{0.2}	\\

\hline


Combined & 50 & & & &
     &      &      & 
     43.6 &      0.0 &
& &
     15.4 &     31.4 &      4.8 &     13.1 &      0.0 
& &
&     15.8 &      8.2 &      1.6\\

no & 100 & & & &

   &      &      & 
    273.7 &     20.6 &
& &
     20.4 &     60.6 &     22.4 &     77.2 &      0.0 & 
& & 
    21.6 &     12.2 &      0.8\\

rack & 150 & & & &

   &     &     &     TL &     35.3 & 

& & 
    26.5 &
     90.4 &     31.4 &    125.2 &      0.6 &
& &
     25.8 &     17.6 &      4.0\\

reduction & 200 & & & &

   &      &      & 
  TL  &     48.7 &
& & 
     32.6 &    122.4 &
     19.8 &    197.1 &      9.9 & 
& &
    37.5 &     17.0 &      1.4\\

&500 & & & &

&      &     & 
    TL &     58.3 &
& & 
     76.5 &    292.4 &      0.2 &    170.7 &     21.0 & 
& &     82.6 &     39.8 &      3.1\\

& \textbf{Average} 	&	
& & &

	& 	& 	&
\textbf{243.4}	& \textbf{32.6}	&
& &
\textbf{34.3}& \textbf{119.4}	& \textbf{15.7}	&

\textbf{116.7}&\textbf{ 6.3}	&
& &
 \textbf{36.7}	& \textbf{19.0}	&\textbf{2.2}

	\\
\hline
Combined	&	50	& & & & & & &	50.4	&	0.0	& & &	15.4	&	31.4	&	4.8	&	12.9	&	0.0	& & &	15.8	&	8.2	&	1.6 \\
rack	&	100	& & & & & & &	252.0	&	19.5	& & &	20.2	&	60.6	&	22.6	&	52.2	&	0.0	& & &	21.7	&	12.2	&	1.0 \\
reduction	&	150	& & & & & & &	TL	&	33.0	& & &	26.2	&	90.4	&	26.6	&	110.7	&	0.9	& & &	26.0	&	17.6	&	4.0 \\
	&	200	& & & & & & &	TL	&	48.4	& & &	33.0	&	122.4	&	21.4	&	143.2	&	9.3	& & &	36.7	&	17.0	&	1.6 \\
	&	500	& & & & & & 	& TL	&	57.9	& & &	77.5	&	292.4	&	3.2	&	168.0	&	19.4	& & &	78.4	&	39.8	&	3.4 \\
	&	\textbf{Average} 	& & & & & & & 	\textbf{240.5}	&	\textbf{31.8}	& & &	\textbf{34.5}	&	\textbf{119.4}	&	\textbf{15.7}	&	\textbf{97.4}	&	\textbf{5.9}	& & &	\textbf{35.7}	&	\textbf{19.0}	&	\textbf{2.3} \\

\hline

\end{tabular}
}
\caption{Order and rack allocation: comparison of average results}
\label{order_compare}
\end{table}

\subsection{Large scale problems}
\label{contrib3}

In terms of RMFS facilities the archetypical example is Amazon, with their large fulfilment centres using Kiva robots. Such facilities might hold upward of hundreds of thousands of different products and contain a correspondingly large number of racks and  robots. It might at first sight appear therefore that the problem to be solved is extremely large.  

However, this is misleading, since for the problem considered in this paper, namely order and rack allocation, the \emph{\textbf{key determinant}} as to the size of the problem to be solved relates not to the size of the fulfilment facility, but rather to \emph{\textbf{the number of unallocated orders ($O$) that have to be allocated in the current pick window}}. 
\cite{allgor2023} give a value of 15 minutes for the pick window in Amazon centres, in other words the pick window only relates to the work that is to be completed within the next 15 minutes.

To investigate how the work given in this paper could be applied to large scale problems, \emph{\textbf{over a single pick window}}, we generated problems involving  up to 40000 products and 10000 racks with 6400 orders to be allocated amongst 64 pickers.
All of these large scale test problems are publicly available at \\
\href{https://www.dcc.ufmg.br/~arbex/mobileRacks.html}{https://www.dcc.ufmg.br/$\sim$arbex/mobileRacks.html}.
Problems of the size considered here are 
significantly larger
than previous problems for multiple pickers considered in the academic research literature. 

For these large scale problems we increased the time limit to 1200 seconds.
\cite{boysen19, boysen19b} report that an average order comprises just 1.6 items. Hence, as a rough approximation, only around $2O$ different products might be involved in $O$ orders.
 Assuming a maximum pick rate for a picker of 600 items per hour (one every 6 seconds,~\cite{lucas}) gives an approximate lower bound on pick time of $(2O \times 6/P)$ seconds for picking $O$ orders using $P$ pickers. This is a lower bound as it assumes that all pickers pick at their maximum rate. For the largest problem in Table~\ref{large} with $O=6400$ and $P=64$ this computes to 1200 seconds. 
Given that (clearly) a RMFS operator can well afford faster hardware than university academics the time limit  we adopted seems reasonable. 

The results for these large scale problems can be seen in Table~\ref{large}.
\cite{valle2021} and Strategy 1,  both of which are  single stage approaches, are not effective on these large scale problems. These single stage approaches fail to produce a feasible solution for three of the eight  problems in Table~\ref{large}. Strategy  2 fails to produce a feasible problem for only one problem, that being the largest problem attempted. Strategy 3  produces  feasible solutions for all problems. For the five problems where the two single stage approaches and Strategy 2 all produce a feasible solution with no backlog, we have that  Strategy 2 requires on average 16.5\% fewer racks than the single stage approaches.

As before Strategy 3, where to streamline rack sequencing  we pick all orders that can be satisfied by a single rack in the first-stage, 
leaving orders that involve two or more racks to the second-stage,  requires more racks than the other approaches. 

With regard to the benefit of rack reduction then, over the large scale problems seen in Table~\ref{large}, we have that our rack reduction procedure reduced the number of racks that need to be considered by an average of 75\%.

\begin{table}[!htpb]
\centering
{\tiny
\renewcommand{\tabcolsep}{1mm} 
\renewcommand{\arraystretch}{1.8} 
\begin{tabular}{|c|c|c|c|c|rrr|rrr|rrrrrrr|rrrrrrr|}
\hline
\multirow[c]{3}{*}{$N$} & \multirow[c]{3}{*}{$R$} & \multirow[c]{3}{*}{$O$} & \multirow[c]{3}{*}{$P$} & \multirow[c]{3}{*}{$C_p$} &
\multicolumn{3}{c|}{\cite{valle2021}} &
 \multicolumn{3}{c|}{Strategy 1} & \multicolumn{7}{c|}{Strategy 2} & \multicolumn{7}{c|}{Strategy 3}\\
\cline{5-25}
& &  &  &  & 
 \multirow[c]{2}{*}{T(s)} & \multirow[c]{2}{*}{GAP} & \multirow[c]{2}{*}{$|\Theta|$} &
\multicolumn{2}{c}{1st stage} & \multirow[c]{2}{*}{$|\Theta|$} & \multicolumn{2}{c}{1st stage} & \multicolumn{2}{c}{2nd stage} & \multirow[c]{2}{*}{$|\Theta|$} & \multirow[c]{2}{*}{$|S|$} & \multirow[c]{2}{*}{BLG} & \multicolumn{2}{c}{1st stage} & \multicolumn{2}{c}{2nd stage} & \multirow[c]{2}{*}{$|\Theta|$} & \multirow[c]{2}{*}{$|S|$} & \multirow[c]{2}{*}{BLG}\\
& &  &  &  &  & & & \multicolumn{1}{c}{T(s)} & \multicolumn{1}{c}{GAP} &  & \multicolumn{1}{c}{T(s)} & \multicolumn{1}{c}{GAP} & \multicolumn{1}{c}{T(s)} & \multicolumn{1}{c}{GAP} &  &  &  & \multicolumn{1}{c}{T(s)} & \multicolumn{1}{c}{GAP} & \multicolumn{1}{c}{T(s)} & \multicolumn{1}{c}{GAP} &  &  & \\
\hline
 \multirow[c]{4}{*}{ 6000} & \multirow[c]{2}{*}{5000}  & \multirow[c]{2}{*}{3200} & 32 & 100 &  
10.4 & -- & 991 &
32.5 &   -- &  991 &  53.0 &   -- &  18.6 & -- &  885 & 1935 &    0 &  24.2 & -- & 67.4 & -- &  507 &  812 &  526 \\
                           &                           &                          & 64 & 50  &  
TL   & 17.2 & 1107 &
TL   & 18.6 & 1108 & 134.9 &   -- & 537.5 & -- & 1002 & 1935 &    0 &  28.8 & -- & 11.0 & -- &  525 &  812 &  643 \\
\cline{2-25}
                           & \multirow[c]{2}{*}{10000} & \multirow[c]{2}{*}{6400} & 32 & 200 &  
12.0 &   -- & 1786 &
33.9 &   -- & 1786 & 156.1 &   -- &   6.9 & -- & 1547 & 3821 &    0 & 151.1 & -- &  6.2 & -- &  940 & 1488 &  962 \\
                           &                           &                          & 64 & 100 &    
TL &      &   NF &
TL &      &   NF & 133.0 &   -- & 127.9 & -- & 1789 & 3821 &    0 & 148.2 & -- & 10.8 & -- & 1037 & 1488 & 1220 \\
\hline
 \multirow[c]{2}{*}{12000} & \multirow[c]{2}{*}{10000} & \multirow[c]{2}{*}{6400} & 32 & 200 & 
39.2 &   -- & 2624 &
136.4 &   -- & 2624 &  80.3 &   -- &  10.5 & -- & 2136 & 3799 &    0 &  78.3 & -- & 18.9 & -- & 1344 & 1802 & 1248 \\
                           &                           &                          & 64 & 100 &    
TL & 23.4 & 2748 &
TL & 24.0 & 2750 &    TL & 0.01 & 382.6 & -- & 2257 & 3799 &    0 &  89.6 & -- & 29.8 & -- & 1389 & 1802 & 1418 \\
 \hline 
 \multirow[c]{2}{*}{40000} & \multirow[c]{2}{*}{10000} & \multirow[c]{2}{*}{6400} & 32 & 200 &    
TL &      &   NF & 
TL &      &   NF & 330.1 &   -- &  58.1 & -- & 2657 & 3805 &    0 &  80.1 & -- & 14.1 & -- & 3054 & 1783 & 1049 \\
                           &                           &                          & 64 & 100 &    
TL &      &   NF &
TL &      &   NF &    TL &      &       &    &   NF & 3805 &      &  75.2 & -- & 19.8 & -- & 3071 & 1783 & 1170 \\
\hline
\end{tabular}
}
\caption{Computational results, rack reduction, large scale problems}
\label{large}
\end{table}


\section{Conclusions}
\label{sec:conclusions}

In this paper we have considered the problem of  allocating orders and mobile storage racks to
multiple pickers 
in a
robotic mobile fulfilment system
when we have inventory constraints and order backlog.

We presented a two stage formulation of the problem. 
In this two stage approach we, in the first-stage, deal with  the orders which must be definitely fulfilled (picked), where the racks chosen to fulfil these 
first-stage orders were chosen so as to (collectively) contain sufficient product to satisfy all orders. In the 
second-stage we restricted attention to the racks chosen in the first-stage solution in terms of allocating second-stage orders.  

We considered three different strategies for first/second-stage order selection. Strategy 1 only used the first-stage 
(so no second-stage) and allocated all orders to pickers. Strategy 2 left the set of orders comprising just a single unit of one product to the second-stage, all other orders were considered in the first-stage. Strategy 3 considered the set of orders which only require one rack as first-stage orders, all other orders were considered in the second-stage. Strategy 3 also ensured that all orders which only require one rack, were all fully supplied by a
single rack, thereby minimising the requirement to make
decisions as to the rack sequence (i.e.~the sequence in which racks are presented to each picker).

We presented a heuristic procedure to reduce the number of racks that need to be considered which was independent of our 
two stage approach and hence could easily be incorporated into any
other solution approach (heuristic or optimal) for order and rack allocation.

Extensive computational results were presented for test problems that are made publicly available;  including test problems that are significantly larger than previous problems considered in the literature.
These results showed that our two stage approach out-performs both a combined single stage approach to order and rack allocation and the approach given previously in~\cite{valle2021}.

Rather than repeat here what we believe to be the contribution to the literature as made by this paper we simply refer the reader back to Section~\ref{sec:contrib}. We believe that the results given in this paper justify the claims made there as to the contribution made by this paper.

 \clearpage
\newpage
 \pagestyle{empty}
\linespread{1}
\small \normalsize 


\appendix
\section{Formulation extensions}
\label{jebappa}
In this Appendix  we briefly consider a number of extensions to the above first and second stage formulations to deal with issues that may arise in practice.

\subsection{Rack distance travelled and storage zones}
In our first-stage formulation the objective adopted addresses minimising the number of racks used. It is a simple matter to change this objective to focus on the rack distance travelled. Let $d_{rp}$ be a measure of the distance travelled if rack $r$ is allocated to picker $p$. Then in Equation~(\ref{eq1}) replace the term relating to the number of racks used,  $\sum_{r=1}^R u_r$, by the total distance travelled, $\sum_{r=1}^R \sum_{p=1}^P d_{rp}y_{rp}$.  

It may be that the  fulfilment centre is  divided into separate zones, where for ease of operation or because of the physical size of the facility,  racks and pickers are dedicated to a single zone, but orders can be allocated to a picker in any zone. If this is the case then this can be easily dealt with by setting $y_{rp}=0$ if rack $r$ cannot be allocated to picker $p$, e.g.~because they are in different storage zones.

\subsection{Picker workload/productivity}
In Equation~(\ref{eq2}) we have defined the capacity of a picker in terms of the number of orders picked, regardless of how long it takes to pick (complete) the order. 
However it is trivial to incorporate picker workload. If $W_o$ is some user-defined measure of the workload associated with order $o$, with (for simplicity of explanation) each picker having the same workload capacity $W^*$ then the total workload for picker $p$ becomes $\sum_{o=1}^O W_o x_{op}$ and we can replace Equation~(\ref{eq2}) in the first-stage formulation by $\sum_{o \in F} W_o x_{op} \leq W^*$. Similarly we replace Equation~(\ref{s2eq2}) in the second-stage formulation by $\sum_{o=1}^O W_o x_{op} \leq W^*$.
Imposing a lower bound on the workload associated with each picker, or having different workload constraints for each picker, in order to maintain picker productivity is also easily done.
Picker workload $W_o$ could, for example, be related to $[q_{io}]$, so the products involved in order $o$ and the number of items of each such product.

\subsection{Rack and order sequencing} 
Equation~(\ref{eq3}) implies that a rack can only be used at most once in the current pick window (if allocated to some picker).
For effective operation of a RMFS it is necessary 
 to have a mechanism (algorithm) for deciding the sequence in which racks are presented to each picker, and too the sequence in which orders are picked from the presented racks. 
If a rack is only used once then it is much easier to achieve a rack and order sequence (since the problem decomposes into deciding $P$ independent sequences). Note here however that,  over time, and multiple pick windows, a rack will typically visit many different pickers. 

In the event that a user wishes to allow racks to be allocated to multiple pickers in the same pick window, the approach given in this paper can be modified appropriately. Introduce a variable (say $\eta_{irp}$) accounting for the number of units of product $i$ stored on rack $r$ that are used at picker $p$. Let $\kappa$ be a limit on the number of pickers to which any individual rack can be allocated. Then in terms of our first-stage formulation,
Equations~(\ref{eq1})-(\ref{eq6}), we remove Equations~(\ref{eq3}),(\ref{eq4}) and add:
\begin{optprog}
& \sum_{p=1}^P y_{rp} & \leq & \kappa & 
\forall r \in \{1,\ldots,R\} \label{jebex4} \\
& \sum_{r=1}^R \eta_{irp} & \geq & \sum_{o \in F}q_{io}x_{op} &
\forall i \in I;~ p \in \{1,\ldots,P\}:~\sum_{o \in F} q_{io} \geq 1
\label{jebex1} \\
& \eta_{irp} & \leq & s_{ir} y_{rp} &
\forall i \in I;~r \in \{1,\ldots,R\}; p \in \{1,\ldots,P\}
\label{jebex11} \\
& \sum_{p=1}^P \eta_{irp} & \leq & s_{ir} & 
\forall i \in I;~r \in \{1,\ldots,R\} \label{jebex2} \\
& y_{rp} & \leq & u_r & 
\forall r \in \{1,\ldots,R\};~p \in \{1,\ldots,P\}\label{jebex3} 
\end{optprog}
Equation~(\ref{jebex4}) limits the number of pickers to which a rack can be allocated to at most 
$\kappa$. This constraint prevents excessive movement of a rack between pickers.
Equation~(\ref{jebex1}) ensures that, for each product $i$ and each picker $p$, sufficient is supplied to meet the required number of units of product $i$ at picker $p$ given the orders allocated to the picker. 
Equation~(\ref{jebex11}) ensures that if a rack is not allocated to a picker (i.e.~$y_{rp}=0$) then no product can be supplied from that rack to the picker (i.e.~we have $\eta_{irp}=0$).
Equation~(\ref{jebex2}) ensures that over all pickers $p$ we do not use more units of product $i$ than are available on rack $r$.
Equation~(\ref{jebex3}) ensures that $u_r$ is one if rack $r$ is used by some picker $p$, and conversely that $y_{rp}$ is zero if $u_r$ is zero.
The corresponding changes needed to the second-stage formulation are very similar to these first-stage changes and so are not given here for space reasons.

\section{Additional constraints}
\label{jebappb}

In this Appendix we detail additional constraints associated with Strategy 3, 
as well as a number of other constraints that we can introduce into the formulation.

\subsection{Additional constraints, Strategy 3}

As discussed above in Strategy 3 
$F$ is the set of orders which can be satisfied (fully supplied) by a single rack (i.e.~the set of orders $o$ with $\delta_o=1$)
and all these orders are picked using just a single rack in the first-stage.
Clearly all orders $o \in Q$ can be fully supplied by a single rack by definition (since any order $o \in Q$ is for a single unit of one product). However to ensure that orders $o \in F \setminus Q$ can be fully supplied using just a single rack we need to add additional constraints to the formulation. 
Let $\Omega(o)$ be the set of racks that can fully supply order $o$, i.e. $\Omega(o) = [r~|~r \in \Delta (o)~ q_{io} \leq s_{ir} ~\forall i \in I]$.
Let $z_{or}=1$ if order $o \in F \setminus Q$ is fully supplied by rack $r \in \Omega(o)$, zero otherwise, where we set
$z_{or}=0$ if $r \notin \Omega(o)$.
 Then we add to the formulation:
\begin{optprog}
& \sum_{r \in \Omega(o)} z_{or} & = & 1  & 
 \forall o \in  F \setminus Q
\label{jebs1} \\
&\sum_{ o \in  F \setminus Q } q_{io}z_{or} & \leq & s_{ir}u_r  
& \forall i \in  I;  r \in \{1,\ldots,R\}
\label{jebs2} \\
& z_{or} - 1 & \leq &  x_{op} - y_{rp} ~ \leq ~ 1 - z_{or}
&  \forall o \in  F \setminus Q;  r \in \Omega(o); p \in \{1,\ldots,P\}
\label{jebs3} \\
& \sum_{r \in \Omega(o)} u_r & \geq & 1 & 
\forall o \in F \setminus Q
\label{jebs4} 
\end{optprog}
Equation~(\ref{jebs1}) ensures that each order $o \in  F \setminus Q$ is allocated to a single rack $ r \in \Omega(o)$. 
 Equation~(\ref{jebs2}) ensures that there is sufficient inventory on a rack to deal with all orders which must be fully supplied from that rack. 
Equation~(\ref{jebs3}) ensures that if $z_{or}=1$ then $ x_{op} = y_{rp}~ p=1,\ldots,P$. In other words if order 
$o$ is fully supplied from rack $r$ then the picker to which order $o$ is (uniquely) allocated corresponds to the picker to which rack $r$ is (uniquely) allocated. If $z_{or}=0$ then this constraint has no effect. Equation~(\ref{jebs4}) ensures that at least one of the racks in $\Omega(o)$ is used for each order  
 $o \in F \setminus Q$. 

We amend  the second-stage formulation for Strategy 3 in a similar manner as for Equations~(\ref{jebs1})-(\ref{jebs3}) discussed above.
Let $z_{or}=1$ if order $o \in F \setminus Q $ is fully supplied by rack $r \in \Omega(o) \cap \Theta $, zero otherwise, where we set
$z_{or}=0$ if $r \notin \Omega(o) \cap \Theta $.  Then we add to the formulation:
\begin{optprog}
& \sum_{r \in \Omega(o) \cap \Theta} z_{or} & = & 1  & 
 \forall o \in  F \setminus Q
\label{jebs21} \\
&\sum_{ o \in  F \setminus Q} q_{io}z_{or} & \leq & s_{ir}  
& \forall i \in  I;  r \in \Theta
\label{jebs22} \\
& z_{or} - 1 & \leq &  x_{op} - y_{rp} ~ \leq ~ 1 - z_{or}
&  \forall o \in  F \setminus Q;  r \in \Omega(o) \cap \Theta; p \in \{1,\ldots,P\}
\label{jebs23} 
\end{optprog}

\subsection{Other constraints}
There are a number of other constraints that we can introduce into the formulation and these are detailed below.

Suppose we have $x_{op}=1$, so order $o$ is allocated to picker $p$. Then it is clear that at least $\delta_o$ of the racks $r \in \Delta(o)$ must be allocated to picker $p$ to fulfil the order. So we have the constraint:
\begin{optprog}
& \sum_{r \in \Delta(o)} y_{rp} & \geq & \delta_o x_{op}   & \forall o \in  \{1,\ldots,O\};  p \in \{1,\ldots,P\}
\label{jeboex1} 
\end{optprog}
Equation~(\ref{jeboex1}) ensures that if order $o$ is allocated to picker $p$ sufficient racks are also allocated to that picker to supply the order. 

Moreover for any product $i$ in order $o$ for which $q_{io} \geq 1$ there must be sufficient  inventory of that product on the racks $r \in \Delta(o)$ allocated to picker $p$. So we have the constraint:
\begin{optprog}
& \sum_{r \in \Delta(o)} s_{ir}y_{rp} & \geq & q_{io} x_{op}  &
& \forall o \in  \{1,\ldots,O\};  p \in \{1,\ldots,P\}; i \in I:~q_{io} \geq 1
\label{jeboex2} 
\end{optprog}

In terms of the first-stage Equations~(\ref{jeboex1}),(\ref{jeboex2}) are amended to restrict attention to orders $o \in F$. In terms of the second-stage 
Equations~(\ref{jeboex1}),(\ref{jeboex2}) are amended to restrict attention to racks $r \in \Delta(o) \cap \Theta$.

It is possible to set a lower bound on the number of orders dealt with by picker $p$ using:
\begin{optprog}
& \sum_{o \in F} x_{op} & \geq & |F| - \sum_{k=1,~k \neq p}^P C_k & \forall p \in \{1,\ldots,P\}
\label{jeblb1} 
\end{optprog}
Equation~(\ref{jeblb1}) follows directly from Equations~(\ref{eq2}),(\ref{eq2a}) where the summation on the right-hand side is the maximum number of orders that can be dealt with excluding picker $p$, leaving at least $|F| - \sum_{k=1,~k \neq p}^P C_k $ orders which must be dealt with using picker $p$. This constraint applies at both stages. Clearly (for picker $p$) it can only be of benefit if $|F| - \sum_{k=1,~k \neq p}^P C_k \geq 1$

We can add some constraints relating to orders which involve one (or more) products that can only be supplied by one rack. Recall that $\Gamma(i)$ is the set of racks that contain product $i$, so $\Gamma(i) = [r~|~s_{ir} \geq 1~ r=1,\ldots,R]$.  In terms of the first-stage if we have some order $o \in F$ containing a product $i$ for which $|\Gamma(i)|=1$ then it is clear that the (unique) rack in $\Gamma(i)$ must be allocated to the same picker as that order. Let $\alpha(i)$ be the first rack in the set $\Gamma(i)$, ties broken arbitrarily. Then we have:
\begin{optprog}

 & x_{op} & = & y_{\alpha(i)p} & \forall o \in F;  p \in \{1,\ldots,P\}:~\exists ~i \in I~s.t.~|\Gamma(i)|=1,~q_{io} \geq 1  
\label{s2new} 
\end{optprog}
Note here that there is no need to set $u_{\alpha(i)}=1$ since that will be enforced by Equation~(\ref{eq4a}).
We would anticipate that at the first-stage Equation~(\ref{s2new}) might be relatively ineffective since in that stage we consider the complete set of racks and so there will be relatively few products $i \in I$ which can only be supplied by one rack. However at the second-stage, where we only consider a subset $\Theta$ of racks, this constraint might be more useful. In addition in the second-stage we can extend Equation~(\ref{s2new})  to orders $o \in S$. 

In more detail, and to ease the notation, suppose that we are considering the second-stage and consequently have set $\Gamma(i) \leftarrow \Gamma(i) \cap \Theta~\forall i \in I$. Then in this second-stage we have that Equation~(\ref{s2new}) applies, and in addition we have 
\begin{optprog}
& v_o - 1 & \leq &  x_{op} - y_{\alpha(i) p} ~\leq ~1 - v_o & 
\forall o \in S;  p \in \{1,\ldots,P\}:~\exists ~i \in I~s.t.~|\Gamma(i)|=1,~q_{io} \geq 1   
\label{s2newa} 
\end{optprog}
Equation~(\ref{s2newa}) is the equivalent of Equation~(\ref{s2new}) for orders $o \in S$ and 
ensures that
if $v_o=1$ then $ x_{op} = y_{\alpha(i)p}~ p=1,\ldots,P$. If $v_o=0$ then this constraint has no effect.

\FloatBarrier


 \clearpage
\newpage
 \pagestyle{empty}
\linespread{1}
\small \normalsize 

\section*{Acknowledgments}
 \noindent Cristiano Arbex Valle was funded by FAPEMIG grant APQ-01267-18.
\noindent Some aspects of this work were undertaken whilst the second author was a Distinguished Research Fellow at NYU London.

\bibliographystyle{plainnat}
\bibliography{tesco}

\end{document}